\documentclass[11pt]{article}

\usepackage[english]{babel}

\usepackage[letterpaper,top=2cm,bottom=2cm,left=3cm,right=3cm,marginparwidth=1.75cm]{geometry}

\usepackage{amsmath} 
\usepackage{amssymb} 
\usepackage{amsthm}   
\usepackage{graphicx}
\usepackage[colorlinks=true, allcolors=blue]{hyperref}
\usepackage{comment}
\usepackage[dvipsnames]{xcolor}
\usepackage{esint}
\usepackage{tikz}
\newtheorem{theorem}{Theorem}[section]
\newtheorem*{theorem*}{Theorem}
\newtheorem{proposition}[theorem]{Proposition}

\newtheorem{lemma}[theorem]{Lemma}

\theoremstyle{remark}
\newtheorem{remark}[theorem]{Remark}

\theoremstyle{definition}
\newtheorem{definition}[theorem]{Definition}

\newcommand{\OmegaT}{\Omega_T}
\newcommand{\uvec}{\mathbf{u}}
\newcommand{\Uvec}{\mathbf{U}}
\newcommand{\psivec}{\mathbf{\Psi}}
\newcommand{\RR}{\mathbb{R}}
\newcommand{\NN}{\mathbb{N}}
\renewcommand{\div}{\mathrm{div}_x}
\newcommand{\nablax}{\nabla_x}
\newcommand{\Svisc}{\mathbb{S}}
\newcommand{\coup}{\alpha}
\newcommand{\paracoup}{\beta}
\newcommand{\Deltax}{\Delta_x}
\newcommand{\dd}{\mathrm{d}}
\newcommand{\Peff}{p_\coup}

\newcommand{\Erel}{\mathcal{E}(\rho,\uvec,c\mid r, \Uvec, C)}

\newcommand{\arho}{\rho_\coup}
\newcommand{\auvec}{\uvec_\coup}
\newcommand{\ac}{c_\coup}
\newcommand{\closOmT}{{\overline{\Omega}_T}}
\newcommand{\curlyE}{\mathcal{E}}
\newcommand{\tgamma}{\tilde{\gamma}}
\newcommand{\rmin}{\underline{r}}
\newcommand{\rmax}{\overline{r}}
\newcommand{\curlyC}{\mathcal{C}}
\newcommand{\cvrate}{s(\coup)}
\newcommand{\cvconst}{K}
\newcommand{\lr}[1]{\left( #1 \right)}
\catcode`\@=11 \@addtoreset{equation}{section}

\catcode`\@=12

\title{Weak-Strong Uniqueness and Relaxation Limit for a Navier--Stokes--Korteweg Model}
\author{Nilasis Chaudhuri$^\dag$, Christian Rohde$ ^* $ and Florian Wendt$ ^* $}

\begin{document}
\date{}
\maketitle

{
\footnotesize
\centerline{$^\dag\;$Institute of Applied Mathematics and Mechanics, University of Warsaw,}
\centerline{ul. Banacha 2 -- 02-097 Warszawa, Poland}

\smallbreak
\centerline{\small \texttt{nchaudhuri@mimuw.edu.pl}}

\bigbreak
\centerline{$^{*}\;$Institute of Applied Analysis and Numerical Simulation, University of Stuttgart,}
\centerline{Pfaffenwaldring 57, 70569 Stuttgart, Germany}
\smallbreak
\centerline{\small \texttt{christian.rohde|florian.wendt@mathematik.uni-stuttgart.de}}
}

\begin{abstract}
We consider a parabolic relaxation model for the compressible Navier--Stokes--Korteweg equations in the isothermal framework.
This system depends on the relaxation parameters $\coup,\paracoup>0$ and approximates formally solutions of the compressible Navier--Stokes--Korteweg equations in the relaxation limit $\coup \to \infty$ and $\paracoup\to 0$.
Introducing the class of finite energy weak solutions for the initial-boundary value problem corresponding to the relaxation model in spatial dimension three, we show that the weak-strong uniqueness principle holds. It asserts that a weak solution and a strong solution emanating from the same initial data coincide as long as the strong solution exists.
Furthermore, we contribute a rigorous convergence result for the relaxation limit $\coup \to \infty$ and $\paracoup\to 0$ and thus justify the relaxation model as an approximate model for the compressible Navier--Stokes--Korteweg equations from a mathematical point of view.
Our results hold for general non-monotone pressure-density relations.
\end{abstract}

\noindent{\textbf{Keywords:} }{Navier-Stokes-Korteweg system; approximation with parabolic relaxation; relative energy; weak-strong uniqueness; relaxation limit.}

\vspace{2mm}

\noindent{\textbf{Mathematics Subject Classification 2020:}}{ 35Q30, 76N06, 35B25, 76T10.  }

\section{Introduction}\label{Sec:Introduction}
\begin{figure}\label{fig:VdW Pressplot}
    \begin{center}
            \begin{tikzpicture}
                \draw[->](0,0) -- (5.6,0) node[right] {$r$};
                \draw[->](0,0) -- (0,5);
                \draw[scale=1.0, domain=0:5.5, smooth, variable=\x, blue] plot ({\x}, {
                0.1*(
                (1.2*\x-2)^3 - 4*(1.2*\x-2)^2 + 24
                )
                });
                \draw[dotted] (1.6667,5.0) -- (1.6667,0) node[below] {$R_1$};
                \draw[dotted] (3.8889,5.0) -- (3.8889,0) node[below] {$R_2$} ;
                \filldraw[color = cyan, opacity=0.2] (0,0) rectangle (1.6667,5.0);
                \filldraw[color = red, opacity=0.2] (3.8889,0.0) rectangle (5.5,5.0);
                \filldraw[color = violet, opacity=0.2] (1.6667,0) rectangle (3.8889, 5.0);
                \draw (0.83,4.4) node[below] {vapor};
                \draw (2.77,4.4) node[below] {spinodal};
                \draw (4.69,4.4) node[below] {liquid};
                \draw (2.77,2.7) node[below] {\color{blue}$p(r)$};
            \end{tikzpicture}
    \end{center}
    \caption{Illustration of a pressure function of Van-der-Waals type.}
\end{figure}
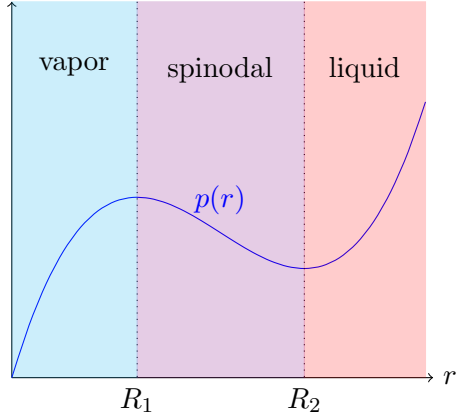

A broadly accepted mathematical description of the dynamics of a compressible viscous single-phase fluid is given by the compressible Navier--Stokes equations.
The extension of this model to a fluid that undergoes phase transition and thus admits two different phases requires to model the phase interfaces appropriately.
In the modeling of these interfaces, there arise two main and substantially different approaches in the literature.
One can either rely on a sharp interface model, where the interface separating two different phases is assumed to be a manifold of codimension one and the fluid's quantities exhibit discontinuities over the interface, or on a diffuse interface model, where the interface is assumed to be a small diffusive region of positive volume and the fluid's quantities vary rapidly but continuously over the interfacial region.
In this paper, we focus on a special instance of a diffuse interface model given by the compressible Navier--Stokes--Korteweg (NSK) equations in an isothermal framework.
We assume that the two-phase fluid occupies a bounded domain $\Omega \subseteq \RR^3$ for some positive time $T>0$.
The NSK equations describe the dynamics of the two-phase fluid by the fluid's density $\rho\colon [0,T] \times \Omega \to \RR_{\geq 0}$ and the fluid's velocity $\uvec\colon [0,T] \times \Omega \to \RR^3$ that obey 
\begin{align}
    &\partial_t \rho 
    +
    \div(\rho\uvec)
    =
    0 && \text{in } (0,T)\times\Omega,\label{NSK field eq continuity}
    \\
    &\partial_t(\rho\uvec)
    +
    \div(\rho\uvec\otimes \uvec) 
    +
    {\nablax}p(\rho)
    -
    \div\Svisc(\nablax \uvec)
    -
    \kappa \rho \nablax \Deltax \rho
    =
    0 &&\text {in } (0,T) \times \Omega.\label{NSK field eq momentum}
\end{align} 
Here, $\kappa>0$ denotes the capillarity coefficient that we assume to be constant in this paper and $p \colon[0,\infty) \to [0,\infty)$ denotes the pressure function.
Furthermore, we denote
\begin{align*}
    \Svisc(A) 
    := 
    2\mu \biggl( \frac{A + A^T}{2} - \frac{1}{3} \mathrm{Tr}(A) \, \mathrm{Id}_{3\times3} \biggr)
    +
    \lambda \mathrm{Tr}(A)\, \mathrm{Id}_{3\times 3} 
\end{align*}
for any $A \in \RR^{3\times 3}$. 
Here, $\mu>0$ and $\lambda>0$ are assumed to be positive constants denoting the shear and bulk viscosity coefficient, respectively.
The term $\div\Svisc(\nablax\uvec)$ in $\eqref{NSK field eq momentum}$ corresponds to assuming that the viscosity mechanism for the two-phase fluid follows Newton's rheological law.
In this paper, we consider zero Neumann boundary conditions for the density and no-slip boundary conditions for the velocity
\begin{align}\label{NSK BC}
    \nablax \rho \cdot \mathbf{n}_{\mid \partial\Omega } = \uvec_{\mid\partial\Omega} = 0 \qquad \text{on } (0,T)\times\partial\Omega
\end{align}
and complement the set of equations $\eqref{NSK field eq continuity}$--$\eqref{NSK BC}$ by initial conditions for $\rho$ and $\uvec$
\begin{align}\label{NSK IC}
    \rho(0) = \rho^0, \qquad 
    \uvec(0) = \uvec^0 \qquad 
    \text{in } \Omega.
\end{align}
The NSK equations are a widely accepted mathematical description of the dynamics of a homogeneous compressible viscous two-phase fluid and their derivation traces back to the seminal work of Korteweg \cite{Korteweg1901}. 
Later, his ideas were picked up by Dunn and Serrin \cite{DunnSerrin1985} and by Anderson, McFadden and Wheeler \cite{AndersonMcFaddenWheeler1998} resulting into the modern formulation $\eqref{NSK field eq continuity}$--$\eqref{NSK field eq momentum}$ of the NSK equations.
The modeling relies on an additional quadratic contribution of the density's gradient in the energy functional that accounts for capillary effects and on a non-convex free energy that results into a non-monotone Van-der-Waals pressure-density relation in order to characterize the two different phases.
The Van-der-Waals pressure-density relation is singular, non-monotone and from a mathematical point on view quite delicate to treat, see e.g.~\cite{FeireislYuNovotny2018}.
Thus, as a mathematical simplification, we only consider pressure functions of Van-der-Waals type in this work, keeping the extension of the subsequent analysis to a Van-der-Waals pressure function open for future work (cf. Section~\ref{Sec:Conclusions}).
Here, we call $p\colon[0,\infty)\to[0,\infty)$ a pressure function of Van-der-Waals type, if there exist two positive constants $0<R_1<R_2<\infty$, such that $p$ is monotonically increasing on $[0,R_1]\cup[R_2,\infty)$ and monotonically decreasing on $(R_1,R_2)$ (cf. Figure~\ref{fig:VdW Pressplot}).
Accordingly, we call then the fluid's state vapor, spinodal and liquid, if the fluid's density value lies in $[0,R_1]$, $(R_1,R_2)$ and $[R_2,\infty)$, respectively.
For more details concerning the modeling we refer to Anderson et al. \cite{AndersonMcFaddenWheeler1998}.\par
The NSK equations has been studied both, from an analytical perspective 
(see e.g.~\cite{HattoriLi1994, BreschDesjardinsLin2003, AntonelliSpirito2022,GermainLeFloch2016})
and from a numerical perspective 
(see e.g.~\cite{DiehlKremserKrönerRohde2016, GiesselmannMakridakisPryer2014,CoquelDiehlMerkleRohde2005}) by many researchers.
In particular, a local-in-time well-posedness result in the class of strong solutions for the initial-boundary value problem has been obtained by Kotschote \cite{Kotschote2008} and the local-in-time existence of classical solutions for the corresponding Cauchy problem has been proved by Hattori and Li \cite{HattoriLi1994}.
Concerning the numerical treatment of the NSK equations with a pressure function of Van-der-Waals type, two major difficulties arise that are not present for the Navier--Stokes equations with a monotone pressure function. On the one hand, an appropriate discretization of a third-order differential operator in the momentum equation is needed.
On the other hand, the underlying first-order system is of elliptic type in the spinodal region so that a straightforward use of classical finite-volume schemes that rely on hyperbolicity is not possible. \par
To overcome these issues, a relaxation system for the NSK equations has been proposed by Rohde \cite{Rohde2010} and was further developed by Hitz et al. \cite{HitzKeimMunzRohde2020}. 
This system approximates the NSK equations by replacing the third-order differential operator by a first-order differential operator via introducing an additional parameter $c \colon[0,T]\times\Omega \to \RR$ that satisfies a linear partial differential equation of second order.
The parameter $c$ is called the order parameter and the precise system reads
\begin{align}
    &\partial_t \rho 
    +
    \div(\rho\uvec)
    = 0 
    &&\text{in } (0,T) \times \Omega,\label{relaxed NSK field eq continuity}
    \\
    &\partial_t(\rho \uvec)
    +
    \div(\rho\uvec\otimes\uvec)
    +
    \nablax p(\rho)
    -
    \div\Svisc(\nablax\uvec)
    -
    \coup \rho  \nablax(c-\rho)
    = 0
    &&\text{in } (0,T) \times \Omega,\label{relaxed NSK field eq momentum}
    \\
    &\paracoup\partial_tc 
    -
    \kappa \Deltax c
    +
    \coup(\rho-c)
    =
    0
    &&\text{in } (0,T)\times\Omega.\label{relaxed NSK field eq parabolic}
\end{align}
Here, the constants $\coup >0$, $\paracoup\geq 0$ are called coupling coefficients.
In the case $\paracoup=0$, the linear equation defining $c$ is elliptic and the resulting system corresponds to the one proposed by Rohde \cite{Rohde2010}.
Based on the elliptic relaxation system $\eqref{relaxed NSK field eq continuity}$--$\eqref{relaxed NSK field eq parabolic}$ with $\paracoup=0$, Hitz et al. \cite{HitzKeimMunzRohde2020} proposed another relaxation system by defining the order parameter $c$ via a parabolic equation in order to prevent a mixed discretization of the governing system of equations.
This system corresponds to the relaxation system for $\paracoup>0$ and was further generalized to a non-isothermal setting by Keim, Munz and Rohde \cite{KeimMunzRohde2023}.
In both cases, one can verify by a formal asymptotic expansion that a sequence of solutions $(\arho,\auvec,\ac)_{\coup>0}$ to the relaxed NSK equations $\eqref{relaxed NSK field eq continuity}$--$\eqref{relaxed NSK field eq parabolic}$ approaches in the relaxation limit $\coup\to\infty$ a solution $(\rho,\uvec)$ of the full NSK equations $\eqref{NSK field eq continuity}$--$\eqref{NSK field eq momentum}$, if $\paracoup=0$ or $\paracoup=\mathcal{O}(\coup^{-1})$ in the case $\paracoup>0$.
For more details concerning the formal asymptotic expansion, we refer to Hitz et al \cite{HitzKeimMunzRohde2020}.
This approximation property was also observed in numerical experiments by Neusser, Rohde and Schleper \cite{NeusserRohdeSchleper2015} for the case $\paracoup=0$ and by Hitz et al. \cite{HitzKeimMunzRohde2020} for the case $\paracoup>0$.
At this point we also refer to Rohde \cite{Rohde2005} for another similar relaxation system that uses a convolutional operator and an interaction potential in order to approximate the third-order differential operator.\par
We complement the relaxed NSK equations $\eqref{relaxed NSK field eq continuity}$--$\eqref{relaxed NSK field eq parabolic}$ by zero Neumann boundary conditions for the order parameter and no-slip boundary conditions for the velocity
\begin{align}\label{relaxed NSK BC}
    \uvec_{\mid\partial\Omega} = \nablax c \cdot \mathbf{n}_{\mid\partial\Omega}=0
    \quad \text{on } (0,T)\times\partial\Omega
\end{align}
and by the initial conditions
\begin{align}\label{relaxed NSK IC}
    \rho(0) = \rho^0,
    \quad 
    \uvec(0)=\uvec^0,
    \quad
    c(0)=c^0
    \quad \text{in }\Omega,
\end{align}
where the initial condition for $c$ is only required in the case $\paracoup>0$.
The momentum equation $\eqref{relaxed NSK field eq momentum}$ can be rewritten as
\begin{align}\label{rewritten momentum relaxed NSK}
    \partial_t(\rho\uvec)
    +
    \div(\rho\uvec\otimes\uvec)
    +
    \nablax \Peff(\rho)
    -
    \div\Svisc(\nablax\uvec)
    -
    \coup \rho \nablax c
    =
    0,
\end{align}
where the artificial pressure function $\Peff$ is defined via
\begin{align*}
    \Peff(r) := p(r) + \frac{\coup}{2} r^2
\end{align*}
for any $r \in [0,\infty)$.
In particular, we have for a pressure function of Van-der-Waals type that the artificial pressure function $\Peff$ is monotonically increasing, provided $\coup>0$ is chosen large enough.
This is an advantage for the numerical treatment of the relaxed NSK equations $\eqref{relaxed NSK field eq continuity}$--$\eqref{relaxed NSK field eq parabolic}$ in comparison to the full NSK equations $\eqref{NSK field eq continuity}$--$\eqref{NSK field eq momentum}$.
\par
From an analytical stand point, only a few results for the relaxed NSK equations $\eqref{relaxed NSK field eq continuity}$--$\eqref{relaxed NSK field eq parabolic}$ are available so far.
In the case $\paracoup=0$, a local-in-time well-posedness result for the Cauchy problem to $\eqref{relaxed NSK field eq continuity}$--$\eqref{relaxed NSK field eq parabolic}$ was obtained by Rohde \cite{Rohde2010}. 
For small initial data, Charve \cite{Charve2014} established a global-in-time well-posedness result and the convergence for the relaxation limit in critical regularity spaces and later extended these results by proving a local-in-time well-posedness result in the same framework \cite{Charve2016}.
Similar results were obtained by Charve and Haspot \cite{CharveHaspot2011, CharveHaspot2013} for the convolutional relaxed NSK equations proposed by Rohde \cite{Rohde2005}.
Furthermore, Haspot \cite{Haspot2010} proved the global-in-time existence of finite energy weak solutions to the corresponding Cauchy problem.
Giesselmann, Lattanzio and Tzavaras \cite[Section 5]{GiesselmannLattanzioTzavaras2017} proved under smoothness assumptions and with periodic boundary conditions a convergence result for the relaxation limit $\coup \to \infty$ in appropriate norms by using a relative energy inequality for the relaxation model proposed by Rohde \cite{Rohde2010}.
\par
It is the goal of this paper to contribute two new analytical results concerning the initial-boundary value problem $\eqref{relaxed NSK field eq continuity}$--$\eqref{relaxed NSK IC}$ for the case $\paracoup>0$.
In our first result, we verify the weak-strong uniqueness property for the initial-boundary value problem $\eqref{relaxed NSK field eq continuity}$--$\eqref{relaxed NSK IC}$ in the class of finite energy weak solutions (cf. Definition~\ref{def:finite energy weak sol}).
More specifically, we show that a finite energy weak solution and a strong solution to $\eqref{relaxed NSK field eq continuity}$--$\eqref{relaxed NSK IC}$ emanating from the same initial data coincide as long as the strong solution exists.
For the compressible Navier--Stokes equations, the weak-strong uniqueness principle for the class of finite energy weak solutions has been obtained by Feireisl, Novotny and Sun \cite{FeireislNovotnySun2011} and by Feireisl, Jin and Novotný \cite{FeireislJinNovotny2012} for a monotone pressure function.
This result was generalized to a compact perturbation of a monotone pressure function by Feireisl \cite{Feireisl2019} and later to a global Lipschitz perturbation by the first author in \cite{Chaudhuri2019}.
Furthermore, Feireisl et al. \cite{FeireislGwiazdaSwierczewskagwiazdaWiedemann2016} established the weak-strong uniqueness principle for dissipative measure-valued solutions with a monotone pressure function, which was later extended to compactly supported perturbations of a monotone pressure function by the first author in \cite{Chaudhuri2020}.
The weak-strong uniqueness principle has various applications in the numerical analysis of the governing equations.
For a detailed overview of these applications concerning the compressible Navier--Stokes equations we refer to \cite{FeireislBook2021}.
With our result we contribute a first step towards obtaining similar applications for the numerical analysis of the relaxed NSK equations $\eqref{relaxed NSK field eq continuity}$--$\eqref{relaxed NSK IC}$.
\par
Our second main result is concerned with the relaxation limit $\coup \to \infty$ in the relaxed NSK equations $\eqref{relaxed NSK field eq continuity}$--$\eqref{relaxed NSK IC}$ for $\beta =\beta(\coup)$ being a function satisfying $\paracoup(\coup) \to 0$ for $\coup \to \infty$.
We prove that a sequence of finite energy weak solutions $(\arho,\auvec,\ac)_{\coup>0}$ to $\eqref{relaxed NSK field eq continuity}$--$\eqref{relaxed NSK IC}$ converges as $\coup \to \infty$ to a strong solution $(\rho,\uvec)$ to the NSK equations $\eqref{NSK field eq continuity}$--$\eqref{NSK IC}$ in certain norms provided that the initial data for $(\arho,\auvec,\ac)$ is not too ill-prepared and as long as the strong solution exists.
As a byproduct, we also verify a convergence rate for the corresponding norms.
In fact, we obtain this convergence result for the limit $\coup \to \infty$ and $\paracoup\to 0$ without assuming any functional relation between $\coup$ and $\paracoup$ (cf. Section~\ref{Sec:Main Results}).
The novelty of this result is two-fold.
On the one hand, we verify this relaxation limit in the class of finite energy weak solutions and without assuming periodicity for the solutions. In particular, we do not assume smoothness for the approximate sequence $(\arho,\auvec,\ac)_{\coup>0}$ generalizing in a certain sense the result in \cite{GiesselmannLattanzioTzavaras2017}.
On the other hand, our result applies for any scaling $\paracoup=\paracoup(\coup)$ that satisfies $\paracoup(\coup) \to 0$ for $\coup \to \infty$, so that our result applies for more general scalings than the ones required by the formal argument via an asymptotic expansion in \cite{HitzKeimMunzRohde2020}.
The relaxation limit can be seen as a singular limit problem for a compressible viscous flow.
We refer to the book by Feireisl and Novotný \cite{FeireislNovotnyBook2009} for various singular limit problems such as the low Mach number transition from compressible to incompressible flows.
\par
As a crucial tool for both our results we derive a relative energy inequality for the relaxed NSK equations $\eqref{relaxed NSK field eq continuity}$--$\eqref{relaxed NSK IC}$.
The leading idea is grounded on the concept of relative entropy tracing back to the pioneering paper of Dafermos \cite{Dafermos1979}.
Adaptions to the compressible Navier--Stokes equations and to the compressible Euler--Korteweg equations can be found in \cite{FeireislJinNovotny2012} and in \cite{GiesselmannLattanzioTzavaras2017}, respectively.
For a relative energy inequality to the compressible NSK equations with density-dependent viscosity coefficients we refer to the recent work by Caggio and Donatelli \cite{CaggioDonatelli2024}.
Analogous to the results for the compressible Navier--Stokes equations, we require a certain amount of monotonicity for the pressure function in order to derive meaningful estimates from the relative energy inequality.
This renders an extension to general pressure functions delicate.
Our results apply for all globally Lipschitz perturbations of a monotone pressure function (cf. $\eqref{assump on p:decomposition}$--$\eqref{assump on p:monotonicity h}$) by adapting the arguments in \cite{Chaudhuri2019}.
In particular, our results account for pressure functions of Van-der-Waals type and thus allow for a two-phase setting.
\\ \par
This paper is organized as follows.
In Section~\ref{Sec:Main Results}, we precisely state the assumptions that we impose on the pressure function throughout this work (cf. $\eqref{assump on p:decomposition}$--$\eqref{assump on p:monotonicity h}$) and define the notion of finite energy weak solutions to the relaxed NSK equations $\eqref{relaxed NSK field eq continuity}$--$\eqref{relaxed NSK IC}$ (cf. Definition~\ref{def:finite energy weak sol}).
Then, we state our two main results given by the weak-strong uniqueness principle (Theorem~\ref{thm main result:W-S uniqueness}) and the relaxation limit (Theorem~\ref{thm main result:singular limit}).
In Section~\ref{Sec:Relative Energy Inequality}, we derive a relative energy inequality for finite energy weak solutions to the relaxed NSK equations $\eqref{relaxed NSK field eq continuity}$--$\eqref{relaxed NSK IC}$ (cf. Proposition~\ref{proposition : relative energy inequality}).
In Section~\ref{Sec:Weak-Strong Uniqueness} we prove our first main result Theorem~\ref{thm main result:W-S uniqueness} and in Section~\ref{Sec:Singular Limit} we prove our second main result Theorem~\ref{thm main result:singular limit} by using the relative energy inequality as a tool.
In Section~\ref{Sec:Conclusions} we close this work with some conclusions.

\section*{Notation}
Let $d \in \NN$ and let $D\subseteq \RR^d$ denote some domain.
For $k \in \NN_0$, we denote the space of $k$-times continuously differentiable functions on $D$ as $\curlyC^k(D)$ and we denote
\begin{align*}
    \|f\|_{\curlyC(D)} := \sup\limits_{x \in D} |f(x)| \quad \text{for } f \in \curlyC^0(D).
\end{align*}
Also, we denote the space of $k$-times continuously compactly supported functions on $D$ as $\curlyC^k_c(D)$ and we set $\curlyC^\infty(D):=\bigcap\limits_{k=0}^\infty \curlyC^k(D)$ and $\curlyC^\infty_c(D) := \bigcap\limits_{k=0}^{\infty} \curlyC^k_c(D)$.
For $p \in [1,\infty]$, we denote by $L^p(D)$ the Lebesgue-spaces on $D$ and by $W^{k,p}(D)$ the $k$-th Sobolev spaces on $D$.
The average of some function $f \in L^1(D)$ is denoted by
\begin{align*}
    \fint_\Omega f \, \dd x 
    := 
    \frac{1}{|D|} \int_D f(x) \, \dd x,
\end{align*}
where $|D|$ denotes the Lebesgue measure on $\RR^d$.
By $\|\cdot\|_{L^p(D)}$ and $\|\cdot\|_{W^{k,p}(D)}$, we denote the Lebesgue-norm and the $k$-th Sobolev norm, respectively.
If $p = 2$, we write $H^k(D)$ instead of $W^{k,2}(D)$.
For $q \in [1,\infty)$, we denote the closure of $\curlyC^\infty_c(D)$ under the norm $\|\cdot\|_{W^{k,q}(D)}$ by $W^{k.q}_0(D)$ and we define $W^{-k,q^\prime}(D)$ as the dual space to $W^{k,q}_0(D)$, where $q^\prime$ denotes the conjugate Hölder exponent to $q$.
For some Banach space $X$ and some positive time $T>0$, we denote the Lebesgue--Bochner spaces on $[0,T]$ ranging into $X$ as $L^p(0,T;X)$ and its norm as $\|\cdot\|_{L^p(0,T;X)}$.
By $\curlyC_{\mathrm{w}}([0,T];L^q(D))$, we denote the space of weakly continuous functions ranging into $L^q(D)$.
For $m \in \NN$, we denote $L^p(D;\RR^m) := L^p(D)^m$ and make the same notation for the other function spaces introduced in this section.
For the space-time cylinder corresponding to $D$, we use the notation
\begin{align*}
    D_T:=(0,T)\times D.
\end{align*}
Finally, for two matrices $A=(A_{ij})_{1\leq i,j\leq 3},\,B = (B_{ij})_{1\leq i,j\leq 3} \in \RR^{3\times 3}$, we denote the Frobenius product as $A:B := \sum\limits_{i,j=1}^{3} A_{ij}B_{ij}$.

\section{Main Results}\label{Sec:Main Results}
This section is devoted to give a precise formulation of our main results.
In order to do so, we first specify the class of pressure functions that we are able to treat in this paper and introduce the class of finite energy weak solutions to the initial-boundary value problem $\eqref{relaxed NSK field eq continuity}$--$\eqref{relaxed NSK IC}$ with $\paracoup>0$.

\subsection{Pressure Function and Finite Energy Weak Solutions}\label{Subsec:Finite Energy Weak Solutions}
We assume that the pressure function $p$ can be decomposed as
\begin{align}\label{assump on p:decomposition}
    p = h + q,
\end{align}
with
\begin{align}\label{assump on p:regularity h and q}
    h \in \curlyC^0([0,\infty))\cap \curlyC^2((0,\infty)),
    \quad 
    q \in \curlyC^{0,1}([0,\infty)) \quad \text{globally Lipschitz},
    \quad 
    q(0) = 0,
\end{align}
and
\begin{align}\label{assump on p:monotonicity h}
    h(0)=0,
    \quad 
    h^\prime>0 \quad \text{on } (0,\infty),
    \quad 
    \lim\limits_{r \to \infty} \frac{h^\prime(r)}{r^{\gamma-1}}=h_\infty
\end{align}
for two positive constants $h_\infty \in (0,\infty)$ and $\gamma \in (1,\infty)$.
These assumptions on the pressure function $p$ correspond essentially to the one considered in \cite{Chaudhuri2019}.
The only difference is that we allow for $\gamma$ the whole range $\gamma \in (1,\infty)$ rather than only $\gamma \in [2,\infty)$.
This is due to the additional quadratic contributions in the energy functional that render the density square integrable, even for $\gamma \in (1,\infty)$ (cf. $\eqref{derivation energy:energy definition at time}$).
The assumptions $\eqref{assump on p:decomposition}$--$\eqref{assump on p:monotonicity h}$ allow for pressure functions of Van-der-Waals type.
Thus, our theory can be applied in the physically relevant case when the relaxation system $\eqref{relaxed NSK field eq continuity}$--$\eqref{relaxed NSK field eq parabolic}$ is used as a model for a compressible viscous liquid-vapor flow.
For a pressure function $p$ satisfying $\eqref{assump on p:decomposition}$--$\eqref{assump on p:monotonicity h}$, we assign pressure potentials $W,H,Q\colon [0,\infty) \to [0,\infty)$ via
\begin{align}\label{pressure potential (def)}
    W(r):= r \int_1^r \frac{p(z)}{z^2} \, \dd z,
    \quad 
    H(r):=r\int_1^r \frac{h(z)}{z^2} \, \dd z,
    \quad 
    Q(r) := r \int_{1}^r \frac{q(z)}{z^2} \, \dd z
\end{align}
for any $r \in [0,\infty)$. 
Furthermore, we readily verify that 
\begin{align}\label{pressure pot:relations I}
    p(r)=W^\prime(r)r - W(r),
    \quad 
    h(r) = H^\prime(r)r - H(r),
    \quad 
    q(r) = Q^\prime(r) r - Q(r),
\end{align}
\begin{align}\label{pressure pot:relations II}
    p^\prime(r) = rW^{\prime\prime}(r),
    \quad 
    h^\prime(r) = rH^{\prime\prime}(r),
    \quad 
    q^\prime(r) = rQ^{\prime\prime}(r)
\end{align}
for any $r \in (0,\infty)$ and
\begin{align}\label{pressure pot:relations III}
    W(r) = H(r) + Q(r)
\end{align}
for any $r \in [0,\infty)$.
Notice that $H$ is a convex pressure potential while both, $W$ and $Q$ are possibly non-convex. 
The decomposition of $W$ into a convex and a non-convex part will be used in Section~\ref{Sec:Relative Energy Inequality}, in order to derive a meaningful relative energy inequality (cf. Lemma~\ref{lemma:auxiliary convexity estimates}).
Due to $\eqref{assump on p:regularity h and q}$ and $\eqref{assump on p:monotonicity h}$ we have that
\begin{align}\label{assump on p:estimates}
    r^\gamma 
    \leq 
    c_1 + c_2 W(r),
    \quad 
    W(r) \leq 
    c_3 (1+r^\gamma)
\end{align}
for any $r \in [0,\infty)$, where $c_1,c_2,c_3>0$ are some positive constants that do not depend on $r$.
In order to introduce the concept of finite energy weak solutions to $\eqref{relaxed NSK field eq continuity}$--$\eqref{relaxed NSK IC}$, we identify the energy for the system and motivate the regularity class for the finite energy weak solution by the following formal calculation.
Let $(\rho,\uvec,c)$ denote some classical solution to the initial-boundary value problem $\eqref{relaxed NSK field eq continuity}$--$\eqref{relaxed NSK IC}$ existing on $[0,T]$.
Multiplying the momentum equation by $\uvec$, integrating over the spatial domain $\Omega$, using the continuity equation $\eqref{relaxed NSK field eq continuity}$ and integration by parts leads to 
\begin{align}\label{derivation energy I}
    \frac{\dd}{\dd t} 
    \int_\Omega 
    \lr{
    \frac{1}{2} \rho |\uvec|^2 
    +
    W(\rho) 
    }
    \, \dd x 
    +
    \int_\Omega
    \Svisc(\nablax\uvec ):\nablax\uvec
    \, \dd x 
    =
    \int_\Omega 
    \coup \rho \uvec \cdot \nablax (c-\rho)
    \, \dd x \, \dd t.
\end{align}
For the right-hand side of this equation, we obtain after using integration by parts, the continuity equation $\eqref{relaxed NSK field eq continuity}$ and the parabolic equation $\eqref{relaxed NSK field eq parabolic}$ 
\begin{align}\label{derivation energy II}
    \int_\Omega 
    \coup \rho\uvec\cdot \nablax (c-\rho)
    \, \dd x
    &=
    \int_\Omega 
    \coup \partial_t \rho \, (c-\rho)
    \, \dd x
    =
    -
    \frac{\dd}{\dd t} 
    \int_\Omega  
    \frac{\coup}{2} |\rho-c|^2 
    \, \dd x
    +
    \int_\Omega 
    \coup \partial_t c \, (c-\rho)
    \, \dd x
\end{align}
and
\begin{align}\label{derivation energy III}
    \int_\Omega 
    \coup 
    \partial_t c \, (c-\rho) 
    \, \dd x
    =
    -
    \int_\Omega 
    \lr{
    \paracoup|\partial_t c|^2 
    -
    \kappa \partial_t c \Deltax c 
    }
    \, \dd x
    =
    -\int_\Omega 
    \paracoup|\partial_t c|^2 
    \, \dd x
    -
    \frac{\dd}{\dd t} 
    \int_\Omega 
    \frac{\kappa}{2}\bigl|\nablax c\bigr|^2 
    \, \dd x.
\end{align}
Inserting $\eqref{derivation energy II}$ and $\eqref{derivation energy III}$ in $\eqref{derivation energy I}$ yields for any $s \in [0,T]$ the energy equality
\begin{align}\label{derivation energy:energy identity}
    \frac{\dd}{\dd t} E(s) 
    +
    \int_\Omega 
    \lr{
    \Svisc(\nablax\uvec) : \nablax \uvec
    +
    \paracoup|\partial_t c|^2 
    }
    \, \dd x 
    = 0,
\end{align}
with 
\begin{align}\label{derivation energy:energy definition at time}
    E(s) 
    := 
    \int_\Omega 
    \lr{
    \frac{1}{2} \rho(s) |\uvec(s)|^2 
    +
    W(\rho(s))
    +
    \frac{\coup}{2} |\rho(s) - c(s)|^2 
    +
    \frac{\kappa}{2} \bigl|\nablax c(s)\bigr|^2 
    }
    \, \dd x.
\end{align}
Integrating for $s \in [0,T]$ the relation $\eqref{derivation energy:energy identity}$ over $(0,s)$ yields
\begin{align}\label{derivation energy:energy identity integrated}
    E(s) 
    +
    \int_0^s \int_\Omega 
    \lr{
    \Svisc(\nablax \uvec):\nablax\uvec
    +
    \paracoup|\partial_t c|^2 
    }
    \, \dd x \, \dd t
    = E_0
\end{align}
with
\begin{align}\label{derivation energy:energy definition initial}
    E_0 
    := 
    \int_\Omega 
    \lr{
    \frac{1}{2} \rho^0|\uvec^0|^2 
    +
    W(\rho^0)
    +
    \frac{\coup}{2} |\rho^0-c^0|^2 
    +
    \frac{\kappa}{2} \bigl|\nablax c^0\bigr|^2 
    }
    \, \dd x.
\end{align}
In what follows let us denote by $\mathcal{K}_0$ a generic positive constant that may vary from line to line but only depends on $\coup,\paracoup,\mu,\lambda,\kappa,T,\Omega$ and on the norms 
\begin{align}\label{derivation energy:initial data norms}
    \Bigl\| \frac{|(\rho\uvec)^0|^2}{\rho^0}\Bigr\|_{L^1(\Omega)},
    \|\rho^0\|_{L^\gamma(\Omega)}, \,
    \|\rho^0-c^0\|_{L^2(\Omega)}, \,
    \|c^0\|_{H^1(\Omega)}.
\end{align}
By $\eqref{assump on p:estimates}$ we have $E_0 \leq \mathcal{K}_0$ and from $\eqref{derivation energy:energy identity integrated}$ we conclude with Hölder's and Poincaré's inequality
\begin{align}
    &\|\sqrt{\rho}\uvec\|_{L^\infty(0,T;L^2(\Omega))} +
    \|\rho\|_{L^\infty(0,T;L^\gamma(\Omega))} +
    \|\rho-c\|_{L^\infty(0,T;L^2(\Omega))} +
    \| \nablax c\|_{L^\infty(0,T;L^2(\Omega))}\nonumber
    \\
    &\qquad 
    +
    \|\uvec\|_{L^2(0,T;H^1(\Omega))}
    +
    \|\partial_t c\|_{L^2(0,T;L^2(\Omega))}
    \leq 
    \mathcal{K}_0.\label{derivation energy:first conclusion from energy identity}
\end{align}
We introduce a function on $[0,T]$ via
\begin{align*}
    \varphi (s) := \int_\Omega c(s) \, \dd x \quad \text{for } s \in [0,T].
\end{align*}
By integrating the parabolic equation $\eqref{relaxed NSK field eq parabolic}$ over the spatial domain over $\Omega$ and using the continuity equation $\eqref{NSK field eq continuity}$, we obtain that the function $\varphi$ satisfies the ODE
\begin{align*}
    &\varphi^\prime(s)
    =
    -\frac{\coup}{\paracoup}\varphi(s)
    +
    \frac{\coup}{\paracoup} 
    \int_\Omega 
    \rho^0
    \, \dd x \quad \text{for } s \in (0,T) ,
    \\
    &\varphi(0)
    = 
    \int_\Omega 
    c^0 
    \, \dd x,
\end{align*}
which has the unique solution
\begin{align*}
    \varphi(s)
    = 
    \Biggl( 
    \int_\Omega 
    c^0 
    \, \dd x 
    +
    \frac{\coup}{\paracoup}\int_\Omega
    \rho^0 \, \dd x
    \int_0^s 
    \exp\Bigl( \frac{\coup}{\paracoup} t \Bigr) 
    \, \dd  t 
    \Biggr)
    \exp\Bigl( -\frac{\coup}{\paracoup} s\Bigr) 
\end{align*}
for any $s \in [0,T]$.
This yields
\begin{align}\label{derivation energy:mean of c}
    \Biggl\|
    \int_\Omega
    c
    \, \dd x
    \Biggr\|_{L^\infty((0,T))}
    \leq 
    \mathcal{K}_0.
\end{align}
From $\eqref{derivation energy:first conclusion from energy identity}$ and $\eqref{derivation energy:mean of c}$ we deduce with Poincaré's inequality
\begin{align}\label{derivation energy:L2 norm of c}
    \|c\|_{L^\infty(0,T;L^2(\Omega))}
    \leq \|\nablax c\|_{L^\infty(0,T;L^2(\Omega))} + \mathcal{K}_0 \Biggl\| \int_\Omega c\, \dd x\Biggr\|_{L^\infty((0,T))}
    \leq \mathcal{K}_0.
\end{align}
Combining $\eqref{derivation energy:first conclusion from energy identity}$ and $\eqref{derivation energy:L2 norm of c}$ yields
\begin{align}\label{derivation energy:L2 estimate density}
    \|\rho\|_{L^\infty(0,T;L^2(\Omega))}
    \leq 
    \|\rho-c\|_{L^\infty(0,T;L^2(\Omega))}
    +
    \|c\|_{L^\infty(0,T;L^2(\Omega))}
    \leq 
    \mathcal{K}_0.
\end{align}
With $\eqref{derivation energy:L2 estimate density}$ we deduce from $\eqref{relaxed NSK field eq parabolic}$ by standard parabolic regularity estimates (see e.g.~\cite[Chapter 6]{Evansbook2010}), that
\begin{align*}
    &\|c\|_{L^2(0,T;H^2(\Omega))}
    +
    \|c\|_{L^\infty(0,T;H^1(\Omega))}
    +
    \|\partial_tc\|_{L^2(0,T;L^2(\Omega))}\\
    &\leq 
    \mathcal{K}_0
    \Bigl( \|c^0\|_{H^1(\Omega)}
    +
    \|\rho\|_{L^2(0,T;L^2(\Omega))} \Bigr)
    \leq \mathcal{K}_0.
\end{align*}
Then, estimates $\eqref{derivation energy:first conclusion from energy identity}$ and $\eqref{derivation energy:L2 estimate density}$ imply for $\tilde{\gamma}:= \max\{ 2,\gamma\}$ by using Hölder's inequality
\begin{align*}
    \|\rho\|_{L^\infty(0,T;L^{\tilde{\gamma}}(\Omega))}
    +
    \|\rho\uvec\|_{L^\infty(0,T;L^{\frac{2\tilde{\gamma}}{\tilde{\gamma}+1}}(\Omega))}
    \leq 
    \mathcal{K}_0.
\end{align*}
By assuming for the initial data $(\rho^0,(\rho\uvec)^0,c^0)$ that the norms in $\eqref{derivation energy:initial data norms}$ are finite, we obtain that the initial data $E_0$ is finite and the preceding estimates motivate the following definition for a finite energy weak solution to the initial-boundary value problem $\eqref{relaxed NSK field eq continuity}$--$\eqref{relaxed NSK IC}$. 

\begin{definition}\label{def:finite energy weak sol}
    Let $\Omega \subseteq \RR^3$ be a bounded domain with regular boundary, $T,\mu,\lambda,\kappa,\coup,\paracoup >0$ and suppose that $p \colon [0,\infty)\to[0,\infty)$ satisfies $\eqref{assump on p:decomposition}$--$\eqref{assump on p:monotonicity h}$ with $\gamma \in (1,\infty)$.
    Let $\tgamma:= \max\{2,\gamma\}$ and let initial data $(\rho^0,(\rho\uvec)^0,c^0)$ with
    \begin{align}
        &\rho^0 \in L^\gamma(\Omega),\quad \rho^0\geq 0 \,\,\, \text{a.e.,} \quad (\rho\uvec)^0 \in L^{\frac{2\tgamma}{\tgamma+1}}(\Omega;\RR^3), \quad c^0 \in H^1(\Omega),\nonumber
        \\
        &\frac{|(\rho\uvec)^0|^2}{\rho^0} \in L^1(\Omega), 
        \quad 
        (\rho\uvec)^0 = 0\,\,\, \text{in } \{\rho^0 = 0\},
        \quad (\rho^0-c^0)\in L^2(\Omega)\label{defi FEW initial conditions}
    \end{align}
    be given.
    Then we call the triplet $(\rho,\uvec,c)$ a finite energy weak solution to $\eqref{relaxed NSK field eq continuity}$--$\eqref{relaxed NSK IC}$ existing on $[0,T]$ emanating from the initial data $(\rho^0,(\rho\uvec)^0,c^0)$, if the following holds:
    \begin{enumerate}
        \item \textbf{Regularity:} we have $\rho\geq 0$ a.e. and
        \begin{align*}
            &\rho \in \curlyC_{\mathrm{w}}([0,T];L^{\tilde{\gamma}}(\Omega)),
            \quad 
            \rho\uvec \in \curlyC_{\mathrm{w}}([0,T];L^{\frac{2\tgamma}{\tgamma+1}}(\Omega;\RR^3)),\\
            &\uvec \in L^2(0,T;H^1_0(\Omega;\RR^3)), 
            \quad 
            c \in L^2(0,T;H^2(\Omega))\cap \curlyC([0,T];H^1(\Omega)),\\
            &\partial_t c \in L^2(0,T;L^2(\Omega)).
        \end{align*}
        \item \textbf{Continuity equation:} We have for any $0\leq s \leq T$ and any test function $\varphi \in \curlyC^1_c([0,T]\times \RR^3)$ 
        \begin{align}
            \biggl[ \int_\Omega \rho(\tau,\cdot)\varphi(\tau,\cdot)\,\dd x \biggr]^{\tau = s}_{\tau=0}
            =
            \int_0^s\int_\Omega 
            \lr{
            \rho \partial_t \varphi 
            +
            \rho \uvec\cdot \nablax\varphi 
            }
            \, \dd x \, \dd t.\label{def:finite energy weak sol:continuity}
        \end{align}
        \item \textbf{Renormalized continuity equation:} We have for any $0 \leq s \leq T$ and any test function $\varphi \in \curlyC^1_c([0,T]\times \RR^3)$ 
        \begin{align}
            &\biggl[ \int_\Omega b(\rho)(\tau,\cdot)  \varphi(\tau,\cdot)\, \dd x  \biggr]^{\tau = s}_{\tau = 0}\nonumber\\
            &\qquad=
            \int_0^s \int_\Omega 
            \lr{
            b(\rho)\partial_t\varphi
            +
            b(\rho) \uvec\cdot \nabla \varphi 
            -
            \bigl( b^\prime(\rho) \rho - b(\rho) \bigr)\div\uvec
            }
            \, \dd x \, \dd t\label{def:finite energy weak sol:renormalized}
        \end{align}
        for any function $b\in \curlyC^1([0,\infty))$ for which there exists some $M_b>0$, such that $b^\prime(z)=0$ for any $z \geq M_b$.
        \item \textbf{Momentum equation:} We have for any $0 \leq s \leq T$ and any test function $\psivec\in \curlyC^1_c([0,T]\times \Omega;\RR^3)$ 
        \begin{align}
            &\biggl[ \int_\Omega \rho\uvec(\tau,\cdot) \psivec(\tau,\cdot) \, \dd x\biggr]^{\tau = s}_{\tau =0}\nonumber\\
            &\qquad =
            \int_0^s \int_\Omega 
            \lr{
            \rho\uvec \cdot \partial_t \psivec
            +
            \rho \uvec\otimes \uvec : \nablax \psivec
            +
            \Peff(\rho)\div\psivec}
            \, \dd x \, \dd t\nonumber\\
            &\qquad \quad -
            \int_0^s \int_\Omega 
            \lr{
            \Svisc(\nablax\uvec): \nablax\psivec
            -
            \coup \rho \nablax c \cdot \psivec
            }
            \, \dd x \, \dd t.\label{def:finite energy weak sol:momentum}
        \end{align}
        \item \textbf{Parabolic equation:} We have for any $0 \leq s \leq T$ and any test function $\varphi \in \curlyC^1_c([0,T]\times \RR^3)$ 
        \begin{align}
            \biggl[ \int_\Omega \paracoup c(\tau,\cdot) \varphi(\tau,\cdot) \, \dd x \biggr]^{\tau = s}_{\tau = 0}
            =
            \int_0^s \int_\Omega 
            \lr{
            \paracoup c \partial_t \varphi 
            -
            \kappa \nabla c \cdot \nabla \varphi 
            -
            \coup(c-\rho)\varphi 
            }
            \, \dd x \, \dd t.\label{def:finite energy weak sol:parabolic}
        \end{align}
        \item \textbf{Energy inequality:} For almost all $s \in (0,T)$, the inequality
        \begin{equation}
            \Bigl[E(s) \Bigr]_{\tau=0}^{\tau=s}
            +
            \int_0^s\int_\Omega 
            \lr{
            \Svisc(\nablax\uvec) : \nablax \uvec
            +
            \paracoup|\partial_t c|^2 
            }
            \, \dd x \, \dd t
            \leq 
            0\label{def:finite energy weak sol:energy inequality}
        \end{equation}
        holds, where $E(s)$ and $E_0$ are defined as in $\eqref{derivation energy:energy definition at time}$ and $\eqref{derivation energy:energy definition initial}$, respectively.
    \end{enumerate}
\end{definition}
\begin{remark}\label{remark:finite energy weak sol}
    \begin{enumerate}
        \item The assumptions $c^0\in H^1(\Omega)$ and $(\rho^0-c^0)\in L^2(\Omega)$ imply $\rho^0 \in L^2(\Omega)$.
        \item The global-in-time existence of finite energy weak solutions for pressure functions $p$ satisfying $\eqref{assump on p:decomposition}$--$\eqref{assump on p:monotonicity h}$ with $\gamma \in (1,\infty)$ can be shown by an adaption of the methods in \cite{AbelsFeireisl2008}.
        Notice that, in contrast to the available results for the compressible Navier--Stokes equations, the global-in-time existence result can be obtained for the full range $\gamma \in (1,\infty)$ due to the quadratic contribution of the density in the energy functional (cf. $\eqref{derivation energy:energy definition at time}$) which renders the density square integrable uniformly-in-time.
        \item Due to the square integrability of the density, the renormalized continuity equation $\eqref{def:finite energy weak sol:renormalized}$ holds in fact for any $b \in \curlyC^0([0,\infty))\cap \curlyC^1((0,\infty))$ that satisfies
        \begin{align*}
            \lim\limits_{r \searrow 0} \Bigl( b^\prime(r)r-b(r) \Bigr) \in \RR,
            \quad 
            |b^\prime(r)| \leq c \quad \forall \, r \in (1,\infty),
        \end{align*}
        for some positive constant $c>0$, see e.g.~\cite[Lemma 11.13]{FeireislNovotnyBook2009}.
        We will need this fact in the derivation of the relative energy inequality (cf. Section~\ref{Sec:Relative Energy Inequality}).
    \end{enumerate}
\end{remark}

\subsection{Weak-Strong Uniqueness}\label{Subsec:Weak-Strong Uniqueness}
Our first main result verifies a weak-strong uniqueness property for finite energy weak solutions to $\eqref{relaxed NSK field eq continuity}$--$\eqref{relaxed NSK IC}$.
\begin{theorem}\label{thm main result:W-S uniqueness}
    Let $\Omega \subseteq \RR^3$ be a bounded domain with regular boundary, $T,\mu,\lambda,\kappa,\coup,\paracoup>0$ and suppose that $p\colon [0,\infty) \to [0,\infty)$ satisfies $\eqref{assump on p:decomposition}$--$\eqref{assump on p:monotonicity h}$ with $\gamma \in (1,\infty)$.
    Let initial data $(\rho^0,(\rho\uvec)^0,c^0)$ satisfying $\eqref{defi FEW initial conditions}$ be given and let $(\rho,\uvec,c)$ denote a finite energy weak solution existing on $[0,T]$ emanating from the initial data $(\rho^0,(\rho\uvec)^0,c^0)$.
    Let the triplet $(r,\Uvec,C)$, with
    \begin{align*}
        &r \in \curlyC^1(\closOmT), \quad 
        \Uvec\in \curlyC^1(\closOmT;\RR^3),
        \quad 
        C \in \curlyC^1(\closOmT),
        \quad r>0,\\
        &\div\Svisc(\nablax\Uvec)\in \curlyC^0(\closOmT;\RR^3),
        \quad 
        \nablax C \in \curlyC^1(\closOmT;\RR^3),
        \quad
        \Deltax C \in \curlyC^1(\closOmT),
    \end{align*}
    be a classical solution to $\eqref{relaxed NSK field eq continuity}$--$\eqref{relaxed NSK IC}$ existing on $[0,T]$ emanating from the same initial data.\\
    Then we have
    \begin{align*}
        \rho(t,x) = r(t,x),
        \quad
        \uvec(t,x) = \Uvec(t,x),
        \quad
        c(t,x) = C(t,x)
        \quad 
    \end{align*}
    for almost all $(t,x) \in (0,T) \times \Omega$.
\end{theorem}
The local-in-time existence of a classical solution to the initial-boundary value problem $\eqref{relaxed NSK field eq continuity}$--$\eqref{relaxed NSK IC}$ is intimately related to appropriate smoothness assumptions imposed on the initial data. 
In the subsequent remark, we discuss local-in-time existence results of strong and classical solutions to the initial-boundary value problem for the compressible Navier--Stokes equations that allow us to anticipate that corresponding results hold also for the relaxed NSK equations.
\begin{remark}
    \begin{enumerate}
    \item 
    The global-in-time existence and uniqueness of classical solutions for the initial-boundary value problem of the compressible Navier--Stokes equations with small initial data has been proved by Matsumura and Nishida \cite{MatsumuraNishida1983}.
    The local-in-time existence and uniqueness of classical solutions for the initial-boundary value problem of the compressible Navier--Stokes equations with initial density away from zero and without any smallness assumption has been obtained by Tani \cite{Tani1977}. By writing the momentum equation of the relaxed NSK equations as in $\eqref{rewritten momentum relaxed NSK}$, we notice that the first two equations $\eqref{relaxed NSK field eq continuity}$--$\eqref{relaxed NSK field eq momentum}$ share the same structure as the compressible Navier--Stokes equations with a force term $\nablax c$.
    By controlling this force term in terms of the density with parabolic regularity estimates, we anticipate that the same method applies for the initial-boundary value problem of the relaxed NSK equations $\eqref{relaxed NSK field eq continuity}$--$\eqref{relaxed NSK IC}$ under the assumption that the initial density is away from vacuum.
    \item 
    Analogous to \cite[Section 3.2.2]{FeireislJinNovotny2012}, it is possible to relax the regularity assumptions on the classical solution $(r,\Uvec,C)$, so that the result is also valid for $(r,\Uvec,C)$ being merely a strong solution of $\eqref{relaxed NSK field eq continuity}$--$\eqref{relaxed NSK IC}$ in appropriate Sobolev spaces.
    The local-in-time existence of a strong solution for the initial-boundary problem of the compressible Navier--Stokes equations with initial density away from zero has been proved by Valli \cite{Valli1983} via a fixed point argument. 
    We anticipate, that the same method applies for the relaxed NSK equations $\eqref{relaxed NSK field eq continuity}$--$\eqref{relaxed NSK IC}$ yielding an analogous local-in-time existence result for strong solutions.
    \end{enumerate}
\end{remark}
We will prove Theorem~\ref{thm main result:W-S uniqueness} in Section~\ref{Sec:Weak-Strong Uniqueness}.
\subsection{Relaxation Limit}\label{Subsec:Model Convergence}
The second main result of this paper tackles the relaxation limit $\coup \to \infty$ of finite energy weak solutions to $\eqref{relaxed NSK field eq continuity}$--$\eqref{relaxed NSK IC}$ with $\paracoup=\paracoup(\coup)$ being a function of $\coup$ satisfying $\paracoup(\coup) \to 0$ for $\coup \to \infty$.
As already discussed in the introduction, we expect that in the limit $\coup \to \infty$, the NSK equations $\eqref{NSK field eq continuity}$--$\eqref{NSK IC}$ will be obtained in a certain sense.
Supposing that the classical solution to the NSK equations $\eqref{NSK field eq continuity}$--$\eqref{NSK IC}$ exists on the same time span as the approximating sequence of finite energy weak solutions, we are able to contribute a rigorous result for this expectation.
The precise statement of our second main result reads as follows.
\begin{theorem}\label{thm main result:singular limit}
    Let $\Omega \subseteq \RR^3$ be a bounded domain with regular boundary and let $T,\mu,\lambda,\kappa,\coup>0$. Suppose that $\paracoup=\paracoup(\coup)$ is a functions of $\coup$ with $\paracoup(\coup) \to 0$ for $\coup \to \infty$.
    Assume that $p\colon [0,\infty) \to [0,\infty)$ satisfies $\eqref{assump on p:decomposition}$--$\eqref{assump on p:monotonicity h}$ with $\gamma \in (1,\infty)$.
    For $\coup>0$, let initial data $(\rho^0_\coup,(\rho\uvec)^0_\coup,c^0_\coup)$ satisfying $\eqref{defi FEW initial conditions}$ be given and let $(\rho_\coup,\uvec_\coup,c_\coup)$ denote a finite energy weak solution of $\eqref{relaxed NSK field eq continuity}$--$\eqref{relaxed NSK IC}$ with $\paracoup=\paracoup(\coup)$ existing on $[0,T]$ emanating from the initial data $(\rho^0_\coup,(\rho\uvec)^0_\coup,c^0_\coup)$.
    Let the couple $(r,\Uvec)$, with
    \begin{align*}
        &r \in \curlyC^1(\closOmT),
        \quad 
        \Uvec\in \curlyC^1(\closOmT;\RR^3),
        \quad 
        r_{\mathrm{min}}:=\inf\limits_{(t,x) \in \closOmT} r(t,x)>0,
        \\
        &\nablax r \in \curlyC^1(\closOmT;\RR^3),
        \quad 
        \Deltax r \in \curlyC^1(\closOmT),
        \quad 
        \div\Svisc(\nablax\Uvec),\,\nablax\div \Uvec\in \curlyC^0(\closOmT;\RR^3),
    \end{align*}
    be a classical solution to $\eqref{NSK field eq continuity}$--$\eqref{NSK IC}$ existing on $[0,T]$ emanating from the initial conditions $(r^0,\Uvec^0)$.
    For $\coup>0$, let
    \begin{align*}
        \mathcal{E}^0_\coup 
        &:=
        \int_\Omega 
        \lr{
        \frac{1}{2} \rho_\coup^0|\uvec_\coup^0 - \Uvec^0|^2
        +
        H(\rho_\coup^0)
        -
        H(r^0)
        -
        H^\prime(r^0) \bigl( \rho_\coup^0 - r^0 \bigr)
        }
        \, \dd x\nonumber
        \\
        &\qquad
        +
        \int_\Omega
        \lr{
        \frac{\coup}{2} |\rho_\coup^0 - c_\coup^0 |^2 
        +
        \frac{\kappa}{2} |\nablax c_\coup^0 - \nablax r^0|^2 
        }
        \, \dd x,
        \\
        \mathcal{E}_\coup 
        &:=
        \int_\Omega 
        \lr{
        \frac{1}{2} \rho_\coup|\uvec_\coup - \Uvec|^2
        +
        H(\rho_\coup)
        -
        H(r^0)
        -
        H^\prime(r) \bigl( \rho_\coup - r \bigr)
        }
        \, \dd x\nonumber
        \\
        &\qquad
        +
        \int_\Omega
        \lr{
        \frac{\coup}{2} |\rho_\coup - c_\coup |^2 
        +
        \frac{\kappa}{2} |\nablax c_\coup - \nablax r|^2 
        }
        \, \dd x.
    \end{align*}
    Then there exists some $\coup_0>0$ such that for any $\coup\geq \coup_0$ the inequality
    \begin{align}\label{thm main result:singular limit:convergence rates}
        &\|\sqrt{\arho}(\auvec-\Uvec)\|_{L^\infty(0,T;L^2(\Omega;\RR^3))}^2
        +
        \|\arho -r \|_{L^\infty(0,T;L^2(\Omega))}^2
        +
        \coup\|\arho-\ac\|_{L^\infty(0,T;L^2(\Omega))}^2\nonumber
        \\
        &\quad +\|\ac-r\|_{L^\infty(0,T;H^1(\Omega))}^2
        +
        \|\auvec-\Uvec\|_{L^2(0,T;H^1(\Omega;\RR^3))}^2\nonumber\\
        &\quad+
        \|\sqrt{\paracoup(\coup)}\partial_t\ac \|_{L^2(0,T;L^2(\Omega))}^2
        \leq 
        \cvconst\cvrate
    \end{align}
    holds with
    \begin{align}\label{definition cvrate and ealpha}
        \cvrate := \frac{1}{\coup} + \paracoup(\coup) + |e_\coup|^2 + \mathcal{E}_\coup^0,\qquad 
        e_\coup:= \fint_\Omega (\arho^0 - r^0) \, \dd x,
    \end{align}
    and where $\cvconst>0$ denotes a positive constant that only depends on $\mu,\lambda,\kappa,\gamma,|\Omega|,T,r_{\mathrm{min}}$ and on the norms
    \begin{align}
        &\|r\|_{\curlyC^1(\closOmT)},
        \|\nablax r\|_{\curlyC^1(\closOmT,\RR^3)},
        \|\Deltax r\|_{\curlyC^1(\closOmT)},\nonumber\\
        &\|\Uvec\|_{\curlyC^1(\closOmT;\RR^3)},
        \|\div\Svisc(\nablax\uvec)\|_{\curlyC^0(\closOmT;\RR^{3})},
        \|\nablax\div\Uvec\|_{\curlyC^0(\closOmT;\RR^3)},\label{thm main result:singular limit:dependencies}
    \end{align}    
    In particular, if 
    \begin{align}\label{thm main result:singular limit:assumption initial data}
        \lim\limits_{\coup \to \infty} \mathcal{E}^0_\coup = 0,
    \end{align}
    then we have the convergences
    \begin{align}
        &\sqrt{\arho}(\auvec-\Uvec) \to 0 \quad \text{in } L^\infty(0,T;L^2(\Omega;\RR^3)),\nonumber
        \\
        &\rho_\coup \to r \quad \text{strongly in } L^\infty(0,T;L^2(\Omega)),\nonumber
        \\
        &c_\coup \to r \quad \text{strongly in } L^\infty(0,T;H^{1}(\Omega)),\nonumber
        \\
        &\sqrt{\coup}(\arho -\ac) \to 0 \quad \text{strongly in } L^\infty(0,T;L^2(\Omega)),\nonumber
        \\
        &\uvec_\coup \to \Uvec \quad \text{strongly in } L^2(0,T;H^{1}(\Omega;\RR^3)),\nonumber
        \\
        &\sqrt{\paracoup(\coup)}\partial_t\ac \to 0 \quad \text{strongly in } L^2(0,T;L^2(\Omega))\label{thm main result:singular limit:convergences}
    \end{align}
    as $\coup \to \infty$.
\end{theorem}
\begin{remark}\label{remark:main result relaxation limit}
    \begin{enumerate}
    \item 
    The local-in-time existence of classical solutions for the Cauchy problem of the NSK equations has been proved in \cite{HattoriLi1994} via the method of successive approximation. 
    We anticipate that the same method also applies for the corresponding initial-boundary value problem $\eqref{NSK field eq continuity}$--$\eqref{NSK IC}$ ensuring the existence of the classical solution in Theorem~\ref{thm main result:singular limit}.
    \item 
    Analogously to \cite[Section 3.2.2]{FeireislJinNovotny2012}, it is possible to relax the regularity assumptions on the classical solution $(r,\Uvec)$ so that the result also holds for $(r,\Uvec)$ being merely a strong solution of $\eqref{NSK field eq continuity}$--$\eqref{NSK IC}$ in appropriate Sobolev spaces.
    The local-in-time existence of a strong solution to $\eqref{NSK field eq continuity}$--$\eqref{NSK IC}$ has been proved by Kotschote \cite{Kotschote2008}. 
    \item Under assumption $\eqref{thm main result:singular limit:assumption initial data}$, inequality $\eqref{thm main result:singular limit:convergence rates}$ yields the convergence rate $\cvrate$ for the corresponding norms on the left-hand side for the limit $\coup \to \infty$.\\
    In particular, if we assume that
    \begin{align*}
        \paracoup(\coup), \mathcal{E}^0_\coup \in \mathcal{O}(\coup^{-1})\quad \text{for } \coup \to \infty,
    \end{align*}
    then we obtain the convergence rate 
    \begin{align*}
        \cvrate\in \mathcal{O}(\coup^{-1}) \quad \text{for } \coup \to \infty.
    \end{align*}
    Note that $\mathcal{E}_\coup^0 \in \mathcal{O}(\coup^{-1})$ implies $|e_\coup|^2 \in \mathcal{O}(\coup^{-1})$ due to the convexity of the potential $H$.
    \item 
    The assumption $\eqref{thm main result:singular limit:assumption initial data}$ and the fact that the pressure potential $H$ is convex imply that
    \begin{align*}
        \rho^0_\coup \to r^0 \quad \text{in } L^{\hat{\gamma}}(\Omega), \quad \text{for } \coup \to \infty,
    \end{align*}
    where $\hat{\gamma}:=\min\{2,\gamma\}$, and
    \begin{align*}
        (\arho^0 - \ac^0) \to 0 \quad \text{in } L^2(\Omega),
        \quad 
        \nablax\ac^0 \to \nablax r^0\quad \text{in } L^2(\Omega;\RR^3) \quad \text{for } \coup \to \infty.
    \end{align*}
    This yields
    \begin{align*}
        \biggl|\fint_\Omega (\ac^0-r^0)\, \dd x \biggr|
        \leq 
        \biggl|\fint_\Omega (\ac^0-\arho^0)\, \dd x \biggr|
        +
        \biggl|\fint_\Omega (\arho^0 - r^0)\, \dd x \biggr|
        \to 0 \quad \text{for } \coup \to \infty
    \end{align*}
    and thus, by Poincaré's inequality,
    \begin{align*}
        \ac^0 \to r^0 \quad \text{in } H^1(\Omega) \quad \text{for } \coup \to \infty.
    \end{align*}
    In particular, we have
    \begin{align*}
        \arho^0 \to r^0 \quad \text{in } L^2(\Omega),\quad
        e_\coup \to 0 \quad \text{for } \coup \to \infty.
    \end{align*}
    \item Inequality $\eqref{thm main result:singular limit:convergence rates}$ and the fact that $\arho$ and $\ac$ satisfy the parabolic equation
    \begin{align*}
        \paracoup(\coup) \partial_t \ac - \kappa \Deltax \ac + \coup (\ac-\arho) = 0
    \end{align*}
    a.e. in $\OmegaT$ imply by Hölder's inequality that
    \begin{align*}
        \Biggl\|\frac{\Deltax \ac}{\sqrt{\coup}} \Biggr\|_{L^2(0,T;L^2(\Omega))}^2 \leq L \cvrate,
    \end{align*}
    where $L>0$ denotes some positive constant that only depends on $\mu,\lambda,\kappa,\gamma,|\Omega|,T,r_{\mathrm{min}}$ and on the norms in $\eqref{thm main result:singular limit:dependencies}$.
    In particular, if $\eqref{thm main result:singular limit:assumption initial data}$ holds, then we have
    \begin{align*}
        \Biggl\|\frac{\Deltax \ac}{\sqrt{\coup}} \Biggr\|_{L^2(0,T;L^2(\Omega))} \to 0 
    \end{align*}
    as $\coup \to \infty$.
    \item 
    In all the numerical experiments from \cite{HitzKeimMunzRohde2020}, the assumptions on the initial conditions were
    \begin{align*}
        \arho^0 = r^0,
        \quad 
        c_\coup^0 = \arho^0,
        \quad 
        \auvec^0 = \Uvec^0.
    \end{align*}
    We readily verify that the assumptions $\eqref{thm main result:singular limit:assumption initial data}$ hold for such a choice of initial data.
    \item The assumption on the functional relation $\paracoup=\paracoup(\coup)$ is not needed. In fact, if we denote for $\coup,\paracoup>0$ a finite energy weak solution $(\rho_{\coup\paracoup},\uvec_{\coup,\paracoup},c_{\coup\paracoup})$ existing on $[0,T]$ emanating from initial conditions $(\rho_{\coup\paracoup}^0,(\rho\uvec)_{\coup\paracoup}^0,c_{\coup\paracoup}^0)$, then there exist $\alpha_0,\paracoup_0>0$, so that an analogous version of inequality $\eqref{thm main result:singular limit:convergence rates}$ holds true for any $\coup\geq \coup_0$ and any $\paracoup\leq \paracoup_0$ (cf. Section~\ref{Sec:Singular Limit}).
    However, owing to the cumbersome notation, we omit such considerations here.
    \end{enumerate}
    We will prove Theorem~\ref{thm main result:singular limit} in Section~\ref{Sec:Singular Limit}.
\end{remark}

\section{The Relative Energy Inequality}\label{Sec:Relative Energy Inequality}
In this section we derive a relative energy inequality for the relaxed system $\eqref{relaxed NSK field eq continuity}$--$\eqref{relaxed NSK IC}$, that is, an inequality that measures a certain distance between regular functions $(r,\Uvec,C)$ satisfying appropriate boundary conditions and a finite energy weak solution $(\rho,\uvec,c)$.
We fix some initial data $(\rho^0,(\rho\uvec)^0,c^0)$ that satisfies $\eqref{defi FEW initial conditions}$ and some finite energy weak solution $(\rho,\uvec,c)$ existing on $[0,T]$ and emanating from this initial data.
Also, we fix some regular functions $r \in \mathcal{C}^1(\closOmT)$, $\Uvec \in \curlyC^1(\closOmT;\RR^3)$ and $C \in \curlyC^1(\closOmT)$ with $r>0$, $\nablax C \in \curlyC^1(\closOmT;\RR^3)$ and $\Deltax C \in \curlyC^1(\closOmT)$ satisfying
\begin{align*}
    \Uvec_{\mid\partial\Omega} = \nablax C\cdot \mathbf{n}_{\mid\partial\Omega} = 0.
\end{align*}
Then, motivated by the results in \cite{FeireislPetcuPravzak2019} and in \cite{Chaudhuri2019}, we consider the relative energy functional
\begin{align*}
    \Erel
    &:= 
    \int_\Omega 
    \lr{
    \frac{1}{2}\rho |\uvec - \Uvec|^2
    +
    H(\rho)
    -
    H(r)
    -
    H^\prime(r) (\rho - r)
    }
    \, \dd x\\
    &\qquad
    +
    \int_\Omega
    \lr{
    \frac{\coup}{2} |(\rho-r) - (c-C)|^2
    +
    \frac{\kappa}{2} |\nablax c-\nablax C|^2
    }
    \, \dd x.
\end{align*}
To obtain a relative energy inequality, we calculate the time evolution of $\Erel$.
To do so, we decompose $\Erel$ by using $\eqref{pressure pot:relations I}$ and $\eqref{pressure pot:relations III}$ as
    \begin{align*}
    \Erel&=\int_\Omega 
    \lr{
        \frac{1}{2} \rho|\uvec|^2 
        +
        W(\rho) 
        +
        \frac{\coup}{2}|\rho - c|^2
        +
        \frac{\kappa}{2}|\nablax c|^2
        }
        \, \dd x\\
        &\quad + 
        \int_\Omega 
        \rho \Bigl( \frac{1}{2}\bigl|\Uvec|^2 - H^\prime(r) \Bigr)\, \dd x
        -\int_\Omega
        (\rho \uvec) \cdot \Uvec
        \, \dd x\\
        &\quad+ \int_\Omega 
        \lr{
        h(r) 
        +
        \frac{\coup}{2}|r-C|^2
        + 
        \frac{\kappa}{2}|\nablax C|^2
        }
        \, \dd x\\
        &\quad -\int_\Omega 
        \kappa \nablax c \cdot \nablax C
        \, \dd x -\int_\Omega 
        \coup (\rho-c)(r-C)
        \, \dd x-\int_\Omega 
        Q(\rho)
        \, \dd x \\ 
        &= 
    \sum\limits_{i=1}^7 I_i
    \end{align*}
    For $I_1$, we use the energy inequality $\eqref{def:finite energy weak sol:energy inequality}$ and obtain for almost any $s \in (0,T)$
    \begin{align*}
        [I_1]^{\tau = s}_{\tau = 0} 
        &\leq 
        -\int_0^s\int_\Omega 
        \lr{
        \Svisc(\nablax\uvec):\nablax \uvec 
        +
        \paracoup |\partial_tc|^2
        }
        \, \dd x \, \dd t.
    \end{align*}
    Thus, we have
    \begin{align*}
        [I_1]_{\tau=0}^{\tau=s}
        +
        \int_0^s \int_\Omega
        \lr{
        \Svisc(\nablax \uvec) : \nablax \uvec
        +
        \paracoup |\partial_t c|^2 
        }
        \, \dd x \, \dd t
        \leq 
        0.
    \end{align*}
    For $I_2$, we test the continuity equation $\eqref{def:finite energy weak sol:continuity}$ with $\frac{1}{2} |\Uvec|^2 - H^\prime(r)$ and obtain
    \begin{align*}
        [I_2]^{\tau=s}_{\tau=0}
        =
        \int_0^s \int_\Omega 
        \rho \Bigl( \Uvec \cdot \partial_t \Uvec - \partial_t H^\prime(r) \Bigr)
        \, \dd x \, \dd t
        +
        \int_0^s \int_\Omega 
        \rho \uvec \cdot \Bigl( \Uvec\cdot \nablax \Uvec - \nablax H^\prime(r) \Bigr)
        \, \dd x \, \dd t.
    \end{align*}
    Thus, we have
    \begin{align*}
        &[I_1 + I_2]_{\tau=0}^{\tau=s} 
        +
        \int_0^s \int_\Omega 
        \lr{
        \Svisc(\nablax \uvec) : \nablax \uvec
        +
        \paracoup |\partial_t c|^2 
        }
        \, \dd x \, \dd t\\
        &\quad \leq 
        \int_0^s 
        \int_\Omega
        \rho\Bigl( \Uvec \cdot \partial_t \Uvec - \partial_t H^\prime(r)  \Bigr)
        \, \dd x \, \dd t
        +
        \int_0^s \int_\Omega
        \rho \uvec \cdot \Bigl( \Uvec \cdot \nablax \Uvec - \nablax H^\prime(r) \Bigr)
        \, \dd x \, \dd t.
    \end{align*}
    For $I_3$ we test the momentum equation $\eqref{def:finite energy weak sol:momentum}$ by $-\Uvec$ and obtain
    \begin{align*}
        [I_3]_{\tau = 0}^{\tau=s}
        &=
        - 
        \int_0^s \int_\Omega 
        \lr{
        \rho \uvec \cdot \partial_t\Uvec 
        +
        \rho \uvec \cdot \uvec \cdot \nabla \Uvec 
        +
        \Bigl( h(\rho) + q(\rho) + \frac{\coup}{2} |\rho|^2 \Bigr) \div \Uvec
        }
        \, \dd x \, \dd t\\
        &\qquad+
        \int_0^s \int_\Omega
        \lr{
        \Svisc(\nablax \uvec) : \nablax \Uvec\, \dd x \, \dd t 
        -
        \coup \rho \nablax c \cdot \Uvec
        }
        \, \dd x \, \dd t.
    \end{align*}
    Thus, we have
    \begin{align*}
        &[I_1 + I_2 + I_3]_{\tau=0}^{\tau=s} 
        +
        \int_0^s \int_\Omega 
        \lr{
        \Svisc(\nablax \uvec) : \bigl(\nablax \uvec- \nablax \Uvec \bigr)
        +
        \paracoup |\partial_t c|^2 
        }
        \, \dd x \, \dd t\\
        &\quad\leq 
        \int_0^s 
        \int_\Omega
        \lr{
        \rho\bigl( \Uvec-\uvec \bigr)  \cdot \partial_t \Uvec 
        +
        \rho \uvec \cdot \bigl( \Uvec - \uvec \bigr) \cdot \nablax \Uvec 
        - 
        \rho\partial_t H^\prime(r)
        -
        \rho\uvec \cdot \nablax  H^\prime(r) 
        }
        \, \dd x \, \dd t\\
        &\qquad
        -
        \int_0^s \int_\Omega 
        \lr{
        \Bigl( h(\rho) + q(\rho) + \frac{\coup}{2} |\rho|^2 \Bigr) \div\Uvec
        +
        \coup \rho \nablax c \cdot \Uvec 
        }
        \, \dd x \, \dd t.
    \end{align*}
    For $I_4$ we rewrite the terms by using the fundamental theorem of integral calculus as
    \begin{align*}
        [I_4]_{\tau=0}^{\tau=s}
        &=
        \int_0^s \int_\Omega 
        \lr{
        \partial_t h(r)
        +
        \frac{\coup}{2}\partial_t(|r-C|^2)
        +
        \frac{\kappa}{2} \partial_t( |\nablax C|^2 )
        }
        \, \dd x \, \dd t\\
        &=
        \int_0^s \int_\Omega 
        \lr{
        h^\prime(r)\partial_t r
        +
        \coup (r-C) \partial_t(r-C)
        -
        \kappa \Deltax C \partial_t C
        }
        \, \dd x \, \dd t,
    \end{align*}
    where we have used integration by parts and the fact that $\nablax C \cdot \mathbf{n}_{\partial\Omega} = 0$ in the last identity.
    Thus, we have
    \begin{align*}
        &[I_1 + I_2 + I_3 + I_4]_{\tau=0}^{\tau=s} 
        +
        \int_0^s \int_\Omega
        \lr{
        \Svisc(\nablax \uvec) : \bigl(\nablax \uvec- \nablax \Uvec \bigr)
        +
        \paracoup |\partial_t c|^2
        }
        \, \dd x \, \dd t\\
        &\quad\leq 
        \int_0^s 
        \int_\Omega 
        \lr{
        \rho\bigl( \Uvec-\uvec \bigr)  \cdot \partial_t \Uvec 
        +
        \rho \uvec \cdot \bigl( \Uvec - \uvec \bigr) \cdot \nablax \Uvec 
        -
        \rho \partial_t H^\prime(r) 
        -
        \rho\uvec \cdot \nablax H^\prime(r) 
        }
        \, \dd x \, \dd t\\
        &\qquad
        -
        \int_0^s \int_\Omega
        \lr{
        \Bigl( h(\rho) + q(\rho) + \frac{\coup}{2} |\rho|^2 \Bigr) \div\Uvec
        +
        \coup \rho \nablax c \cdot \Uvec
        }
        \, \dd x \, \dd t
        \\
        &\qquad
        +
        \int_0^s \int_\Omega 
        \lr{
        h^\prime(r) \partial_t r
        +
        \coup (r-C)\partial_t(r-C)
        -
        \kappa \Deltax C \partial_t C
        }
        \, \dd x \, \dd t.
    \end{align*}
    Concerning $I_5$, we proceed analogously and obtain
    \begin{align*}
        [I_5]_{\tau=0}^{\tau=s}
        &=
        \Biggl[ \int_\Omega \kappa c\Deltax C \, \dd x \Biggr]_{\tau=0}^{\tau=s}
        =
        \int_0^s \int_\Omega 
        \kappa \partial_t(  c \Deltax C) 
        \, \dd x \, \dd t\\
        &=
        \int_0^s \int_\Omega 
        \lr{
        \kappa \partial_t c\Deltax C
        +
        \kappa c \partial_t \Deltax C
        }
        \, \dd x \, \dd t
        =
        \int_0^s \int_\Omega 
        \lr{
        \kappa \partial_t c \Deltax C 
        + 
        \kappa \Deltax c \partial_t C 
        }
        \, \dd x \, \dd t,
    \end{align*}
    where we have integration by parts and the fact that $\nablax C\cdot\mathbf{n}_{\partial\Omega} = 0$ in the first a last identity.
    Thus, we have
    \begin{align*}
        &[I_1 + I_2 + I_3 + I_4 + I_5]_{\tau=0}^{\tau=s} 
        +
        \int_0^s \int_\Omega
        \lr{
        \Svisc(\nablax \uvec) : \bigl(\nablax \uvec- \nablax \Uvec \bigr)
        +
        \paracoup |\partial_t c|^2
        }
        \, \dd x \, \dd t\\
        &\leq 
        \int_0^s 
        \int_\Omega 
        \lr{
        \rho\bigl( \Uvec-\uvec \bigr)  \cdot \partial_t \Uvec 
        +
        \rho \uvec \cdot \bigl( \Uvec - \uvec \bigr) \cdot \nablax \Uvec 
        -
        \rho \partial_t  H^\prime(r) 
        -
        \rho\uvec\cdot \nablax H^\prime(r) 
        }
        \, \dd x \, \dd t\\
        &\qquad
        -
        \int_0^s \int_\Omega
        \lr{
        \Bigl( h(\rho) + q(\rho) + \frac{\coup}{2} |\rho|^2 \Bigr) \div\Uvec
        -
        \coup \rho \nablax c \cdot \Uvec
        }
        \, \dd x \, \dd t
        \\
        &\qquad
        +
        \int_0^s \int_\Omega 
        \lr{
        h^\prime(r) \partial_t r
        +
        \coup(r-C)\partial_t(r-C)
        -
        \kappa\Deltax C \partial_t C
        +
        \kappa \partial_t c \Deltax C 
        +
        \kappa \Deltax c \partial_t C
        }
        \, \dd x \, \dd t.
    \end{align*}
    Concerning $I_6$, we first decompose
    \begin{align*}
        [I_6]_{\tau=0}^{\tau=s}
        &=
        \biggl[ -\int_\Omega \coup \rho(r-C)\, \dd x \biggr]_{\tau=0}^{\tau =s}
        +
        \biggl[  \int_\Omega \coup c(r-C) \, \dd x \biggr]^{\tau=s}_{\tau=0}.
    \end{align*}
    For the first term, we test the continuity equation $\eqref{def:finite energy weak sol:continuity}$ by $-\coup(r-C)$ and obtain
    \begin{align*}
        \biggl[ -\int_\Omega \coup \rho(r-C)\, \dd x \biggr]_{\tau=0}^{\tau =s}
        &=
        -\int_0^s\int_\Omega
        \lr{
        \coup \rho \partial_t(r-C) 
        +
        \coup \rho \uvec \cdot \nablax( r-C)
        }
        \, \dd x \, \dd t.
    \end{align*}
    For the second term, we use again the fundamental theorem of integral calculus to obtain
    \begin{align*}
        \biggl[  \int_\Omega \coup c (r-C) \, \dd x  \biggr]_{\tau=0}^{\tau=s}
        &=
        \int_0^s \int_\Omega 
        \coup \partial_t \bigl( c(r-C) \bigr)
        \, \dd x \, \dd t\\
        &=
        \int_0^s \int_\Omega 
        \lr{
        \coup\partial_tc\,(r-C) 
        +
        \coup c \partial_t( r-C)
        }
        \, \dd x \, \dd t.
    \end{align*}
    Together, this yields
    \begin{align*}
        [I_6]_{\tau = 0}^{\tau = s}
        &=
        \int_0^s \int_\Omega 
        \lr{
        \coup(c-\rho)\partial_t(r-C)
        -
        \coup \rho \nablax (r-C) \cdot \uvec
        +
        \coup (r-C) \partial_t c
        }
        \, \dd x \, \dd t.
    \end{align*}
    Thus, we have
    \begin{align*}
        &[I_1 + I_2 + I_3 + I_4 + I_5 + I_6]^{\tau=s}_{\tau=0}
        +
        \int_0^s \int_\Omega 
        \lr{
        \Svisc(\nablax \uvec) : \bigl(\nablax \uvec- \nablax \Uvec \bigr)
        +
        \paracoup |\partial_t c|^2 
        }
        \, \dd x \, \dd t\\
        &\leq 
        \int_0^s 
        \int_\Omega 
        \lr{
        \rho\bigl( \Uvec-\uvec \bigr)  \cdot \partial_t \Uvec 
        + 
        \rho \uvec \cdot \bigl( \Uvec - \uvec \bigr) \cdot \nablax \Uvec 
        -
        \rho \partial_t H^\prime(r) 
        -
        \rho\uvec \cdot \nablax  H^\prime(r)
        }
        \, \dd x \, \dd t\\
        &\qquad
        -
        \int_0^s \int_\Omega 
        \lr{
        \Bigl( h(\rho) + q(\rho) + \frac{\coup}{2}|\rho|^2 \Bigr) \div\Uvec
        +
        \coup \rho \nablax c \cdot \Uvec
        }
        \, \dd x \, \dd t
        \\
        &\qquad
        +
        \int_0^s \int_\Omega 
        \lr{
        h^\prime(r) \partial_t r
        +
        \coup (r-C) \partial_t (r-C) 
        -
        \kappa \Deltax C \partial_t C
        +
        \kappa \partial_t c \Deltax C 
        +
        \kappa \Deltax c \partial_t C
        }
        \, \dd x \, \dd t\\
        &\qquad+
        \int_0^s \int_\Omega 
        \lr{
        \coup (c-\rho)\partial_t(r-C)
        -
        \coup \rho \nablax (r-C) \cdot \uvec
        +
        \coup(r-C) \partial_t c
        }
        \, \dd x \, \dd t.
    \end{align*}
    For $I_7$, we first verify by using the Lipschitz continuity of $q$ and $\eqref{pressure pot:relations I}$ that $Q\in \curlyC^0([0,\infty))\cap \curlyC^1((0,\infty))$ and that
    \begin{align*}
        \lim\limits_{r \searrow 0} \Bigl(Q^\prime(r)r - Q(r) \Bigr) 
        =
        \lim\limits_{r \searrow 0} q(r) = q(0)\in \RR,
        \quad 
        |Q^\prime(r)| \leq K
        \quad \forall \, r \in (1,\infty)
    \end{align*}
    for some positive constant $K>0$ that does not depending on $r$.
    Thus, according to Remark~\ref{remark:finite energy weak sol}, the renormalized continuity equation $\eqref{def:finite energy weak sol:renormalized}$ holds for $b=Q$. 
    By specifically choosing $\varphi \equiv 1$ as a test function, we obtain 
    \begin{align*}
        \bigl[ I_7 \bigr]^{\tau=s}_{\tau=0}
        =
        -\Biggl[ \int_\Omega Q(\rho) \, \dd x  \Biggr]_{\tau=0}^{\tau=s}
        =
        \int_0^s \int_\Omega 
        q(\rho) \div \uvec 
        \, \dd x \, \dd t.
    \end{align*}
        Thus, we have
    \begin{align*}
        &[\Erel]^{\tau=s}_{\tau=0}
        +
        \int_0^s \int_\Omega
        \lr{
        \Svisc(\nablax \uvec - \nablax \Uvec) : \bigl(\nablax \uvec- \nablax \Uvec \bigr)
        +
        \paracoup |\partial_t c|^2 
        }
        \, \dd x \, \dd t
        \\
        &\leq 
        \int_0^s \int_\Omega 
        \lr{
        \Svisc(\nablax \Uvec ) : (\nablax \Uvec - \nablax \uvec )
        }
        \, \dd x \, \dd t
        \\
        &\qquad+
        \int_0^s 
        \int_\Omega 
        \lr{
        \rho\bigl( \Uvec-\uvec \bigr)  \cdot \partial_t \Uvec 
        + 
        \rho \uvec \cdot \bigl( \Uvec - \uvec \bigr) \cdot \nablax \Uvec 
        -
        \rho \partial_t H^\prime(r) 
        -
        \rho\uvec \cdot \nablax H^\prime(r) 
        }
        \, \dd x \, \dd t\\
        &\qquad
        -
        \int_0^s \int_\Omega 
        \lr{
        \Bigl( h(\rho) + \frac{\coup}{2}|\rho|^2 \Bigr) \div\Uvec
        +
        q(\rho) \Bigl( \div\Uvec - \div\uvec \Bigr)
        +
        \coup \rho \nablax c \cdot \Uvec
        }
        \, \dd x \, \dd t
        \\
        &\qquad
        +
        \int_0^s \int_\Omega 
        \lr{
        h^\prime(r) \partial_t r
        +
        \coup (r-C) \partial_t (r-C) 
        -
        \kappa \Deltax C \partial_t C
        +
        \kappa \partial_t c \Deltax C 
        +
        \kappa \Deltax c \partial_t C
        }
        \, \dd x \, \dd t\\
        &\qquad+
        \int_0^s \int_\Omega 
        \lr{
        \coup (c-\rho)\partial_t(r-C)
        -
        \coup \rho \nablax (r-C) \cdot \uvec
        +
        \coup(r-C) \partial_t c
        }
        \, \dd x \, \dd t.
    \end{align*}
    Introducing the remainders 
    \begin{align*}
        \mathcal{R}^1:&=
        \int_0^s\int_\Omega 
        \biggl(
        \rho \Bigl( \partial_t \Uvec + \uvec\cdot \nablax \Uvec \Bigr) \cdot \Bigl( \Uvec - \uvec \Bigr)
        +
        \Svisc(\nablax \Uvec):\Bigl( \nablax\Uvec - \nablax \uvec \Bigr)
        \biggr)
        \, \dd x \, \dd t
        \\
        &\qquad +
        \int_0^s \int_\Omega 
        \lr{
        h^\prime(r) \partial_t r
        -
        h(\rho) \div\Uvec 
        -
        \rho \partial_t H^\prime(r)
        -
        \rho \uvec \nablax H^\prime(r)
        }
        \, \dd x \, \dd t
        \\
        &\qquad 
        +
        \int_0^s \int_\Omega 
        q(\rho) \Bigl( \div\uvec - \div\Uvec \Bigr)
        \, \dd x \, \dd t,
        \\
        \mathcal{R}^2&:=
        -
        \int_0^s \int_\Omega 
        \lr{
        \frac{\coup}{2}|\rho|^2 \div \Uvec
        +
        \coup \rho \nablax c \cdot \Uvec
        +
        \coup\rho\nablax(r-C)\cdot\uvec
        }
        \, \dd x \, \dd t,
        \\
        \mathcal{R}^3&:=
        \int_0^s \int_\Omega 
        \coup \Bigl( (c-\rho) - (C-r)  \Bigr) \partial_t r
        \, \dd x \, \dd t,
        \\
        \mathcal{R}^4&:=
        \int_0^s \int_\Omega
        \biggl(
        \Bigl( \kappa \Deltax c - \coup (c-\rho) \Bigr) \partial_t C
        -
        \Bigl( \kappa \Deltax C - \coup (C - r) \Bigr) \partial_t C
        \biggr) \, \dd x \, \dd t\nonumber\\
        &\qquad+
        \int_0^s \int_\Omega 
        \Bigl( \kappa \Deltax C - \coup (C-r) \Bigr) \partial_t c
        \, \dd x \, \dd t,
    \end{align*}
    we concisely write
    \begin{align*}
        &[\Erel]^{\tau = s}_{\tau = 0}
        +
        \int_0^s \int_\Omega 
        \lr{
        \Svisc(\nablax \uvec - \nablax \Uvec) : (\nablax \uvec - \nablax \Uvec)
        +
        \paracoup |\partial_t c|^2
        }
        \, \dd x \, \dd t\\
        &\qquad\leq 
        \mathcal{R}^1 
        +
        \mathcal{R}^2
        +
        \mathcal{R}^3
        +
        \mathcal{R}^4.
    \end{align*}
    Let us now rewrite the classical remainder $\mathcal{R}^1$ in such a way that we can see where the momentum equation for a smooth solution will appear later on.
    We first deduce from $\eqref{pressure pot:relations II}$ 
    \begin{equation*}
        h^\prime(r) \partial_t r
        =
        H^{\prime \prime}(r) r \partial_t r
        =
        r\partial_t H^\prime(r)
    \end{equation*}
    and
    \begin{equation*}
        \nablax h(r)
        =
        h^\prime(r) \nablax r
        =
        H^{\prime\prime}(r) r \nablax r
        =
        r \nablax H^\prime(r).
    \end{equation*}
    With these identities, we conclude
    \begin{equation}\label{Rewriting R^1 I}
        \int_0^s \int_\Omega 
        \lr{
        h^\prime(r)\partial_t r - \rho \partial_t H^\prime(r)
        }
        \, \dd x \, \dd t
        =
        \int_0^s \int_\Omega 
        \lr{
        (r-\rho) \partial_t H^\prime(r)
        }
        \, \dd x \, \dd t
    \end{equation}
    and
    \begin{align}
        -\int_0^s\int_\Omega 
        \rho  \uvec \cdot \nablax H^\prime(r)
        \, \dd x \, \dd t
        &=
        \int_0^s \int_\Omega 
        \biggl(
        \Bigl( r \Uvec-\rho \uvec\Bigr) \cdot \nablax H^\prime(r)
        -
        r\Uvec\cdot \nablax H^\prime(r) 
        \biggr)
        \, \dd x \, \dd t\nonumber\\
        &=
        \int_0^s \int_\Omega 
        \biggl(
        \Bigl( r \Uvec - \rho  \uvec \Bigr) \cdot \nablax H^\prime(r)
        -
        \Uvec \cdot \nablax h(r)
        \biggr)
        \, \dd x \, \dd t\nonumber\\
        &=
        \int_0^s \int_\Omega
        \biggl(
        \Bigl( r \Uvec - \rho \uvec \Bigr) \cdot \nablax H^\prime(r)
        +
        h(r) \div\Uvec
        \biggr)
        \, \dd x \, \dd t,\label{Rewriting R^1 II}
    \end{align}
    where we have used integration by parts and the fact that $\Uvec_{\mid\partial\Omega}=0$ to obtain the last identity.
    Finally, we rewrite the convective term as
    \begin{align}
        &\int_0^s \int_\Omega 
        \rho\Bigl( \partial_t \Uvec + \uvec\cdot \nablax \Uvec \Bigr) \cdot \Bigl( \Uvec-\uvec \Bigr)
        \, \dd x \, \dd t\nonumber
        \\
        &\quad 
        =
        \int_0^s \int_\Omega 
        \rho \Bigl( \partial_t \Uvec + \Uvec \cdot \nablax \Uvec \Bigr) \cdot \Bigl( \Uvec-\uvec \Bigr)
        +
        \rho  \Bigl(\uvec - \Uvec \Bigr) \cdot \nablax \Uvec \cdot \Bigl( \Uvec - \uvec\Bigr) 
        \, \dd x \, \dd t.\label{Rewriting R^1 III}
    \end{align}
    Combining $\eqref{Rewriting R^1 I}, \eqref{Rewriting R^1 II}$ and $\eqref{Rewriting R^1 III}$ yields
    \begin{align*}
        \mathcal{R}^1 &=
        \int_0^s \int_\Omega 
        \biggl(
        \rho \Bigl( \partial_t \Uvec + \Uvec \cdot \nablax \Uvec \Bigr) \Bigl( \Uvec - \uvec \Bigr)
        +
        \rho \Bigl( \uvec - \Uvec \Bigr) \cdot \nablax \Uvec \cdot \Bigl( \Uvec - \uvec \Bigr)
        \biggr)
        \, \dd x \, \dd t\nonumber\\
        &\qquad +
        \int_0^s \int_\Omega 
        \biggl(
        \Svisc(\nablax \Uvec) : \Bigl( \nablax \Uvec - \nablax \uvec \Bigr)
        +
        q(\rho) \Bigl( \div\uvec - \div \Uvec \Bigr)
        \biggr)
        \, \dd x \, \dd t\nonumber\\
        &\qquad +
        \int_0^s \int_\Omega 
        \biggl(
        \Bigl( h(r) - h(\rho) \Bigr) \div\Uvec 
        +
        (r-\rho) \partial_t H^\prime(r)
        +
        (r\Uvec - \rho \uvec) \cdot \nablax H^\prime(r)
        \biggr)
        \, \dd x \, \dd t.
    \end{align*}
    To sum up, we have shown the following relative energy inequality that will be our tool to prove our two main results Theorem~\ref{thm main result:W-S uniqueness} and Theorem~\ref{thm main result:singular limit}.
    \begin{proposition}[Relative energy inequality]\label{proposition : relative energy inequality}
        Let initial conditions $(\rho^0,(\rho\uvec)^0,c^0)$ be given satisfying $\eqref{defi FEW initial conditions}$ and let $(\rho,\uvec,c)$ be a finite energy weak solution existing on $[0,T]$ emanating from the initial conditions $(\rho^0,(\rho\uvec)^0,c^0)$.
        Let regular functions $r \in \curlyC^1(\closOmT)$, $\Uvec \in \curlyC^1(\closOmT;\RR^3)$ and $C \in \curlyC^1(\closOmT)$ with $r>0$, $\nablax C \in \curlyC^1(\closOmT;\RR^3)$ and $\Deltax C \in \curlyC^1(\closOmT)$ be given satisfying
        \begin{align*}
            \Uvec_{\mid\partial\Omega} = \nablax C \cdot \mathbf{n}_{\mid\partial\Omega} = 0.
        \end{align*}
        Then the relative energy inequality
        \begin{align}
            &[\Erel]^{\tau = s}_{\tau = 0}
            +
            \int_0^s \int_\Omega 
            \lr{
            \Svisc(\nablax \uvec - \nablax \Uvec) : (\nablax \uvec - \nablax \Uvec)
            +
            \paracoup |\partial_t c|^2
            }
            \, \dd x \, \dd t\nonumber\\
            &\qquad \leq 
            \mathcal{R}^1 + \mathcal{R}^2 + \mathcal{R}^3 + \mathcal{R}^4\label{proposition : relative energy inequality : inequality}
        \end{align}
        holds for almost all $s \in (0,T)$, where
        \begin{align*}
        \mathcal{R}^1 &:=
        \int_0^s \int_\Omega 
        \biggl(
        \rho \Bigl( \partial_t \Uvec + \Uvec \cdot \nablax \Uvec \Bigr) \Bigl( \Uvec - \uvec \Bigr)
        +
        \rho \Bigl( \uvec - \Uvec \Bigr) \cdot \nablax \Uvec \cdot \Bigl( \Uvec - \uvec \Bigr)
        \biggr)
        \, \dd x \, \dd t\nonumber\\
        &\qquad +
        \int_0^s \int_\Omega 
        \biggl(
        \Svisc(\nablax \Uvec) : \Bigl( \nablax \Uvec - \nablax \uvec \Bigr)
        +
        q(\rho) \Bigl( \div\uvec - \div \Uvec \Bigr)
        \biggr)
        \, \dd x \, \dd t\nonumber\\
        &\qquad +
        \int_0^s \int_\Omega 
        \biggl(
        \Bigl( h(r) - h(\rho) \Bigr) \div\Uvec 
        +
        (r-\rho) \partial_t H^\prime(r)
        +
        (r\Uvec - \rho \uvec) \cdot \nablax H^\prime(r)
        \biggr)
        \, \dd x \, \dd t,
        \\
        \mathcal{R}^2&:=
        -\int_0^s \int_\Omega 
        \lr{
        \frac{\coup}{2}|\rho|^2 \div \Uvec
        +
        \coup \rho \nablax c \cdot \Uvec
        +
        \coup\rho\nablax(r-C)\cdot\uvec
        }
        \, \dd x \, \dd t,
        \\
        \mathcal{R}^3&:=
        \int_0^s \int_\Omega 
        \coup \Bigl( (c-\rho) - (C-r)  \Bigr) \partial_t r
        \, \dd x \, \dd t,
        \\
        \mathcal{R}^4&:=
        \int_0^s \int_\Omega 
        \biggl(
        \Bigl( \kappa \Deltax c - \coup (c-\rho) \Bigr) \partial_t C
        -
        \Bigl( \kappa \Deltax C - \coup (C - r) \Bigr) \partial_t C
        \biggr)
        \, \dd x \, \dd t\nonumber
        \\
        &\qquad
        +
        \int_0^s \int_\Omega
        \Bigl( \kappa \Deltax C - \coup (C-r) \Bigr) \partial_t c
        \, \dd x \, \dd t.
        \end{align*}
    \end{proposition}
    
    \section{Weak-Strong Uniqueness}\label{Sec:Weak-Strong Uniqueness}
    The goal of this section is to provide the proof of our first main result Theorem~\ref{thm main result:W-S uniqueness}.
    As a tool, we will use the relative energy inequality derived in Section~\ref{Sec:Relative Energy Inequality}.
    More specifically, we fix initial data $(\rho^0,(\rho\uvec)^0,c^0)$ satisfying $\eqref{defi FEW initial conditions}$ and a finite energy weak solution $(\rho,\uvec,c)$ of $\eqref{relaxed NSK field eq continuity}$--$\eqref{relaxed NSK IC}$ existing on $[0,T]$ emanating from the initial data $(\rho^0,(\rho\uvec)^0,c^0)$.
    Also, we assume that functions $r \in \curlyC^1(\closOmT)$ with $r>0$, $\Uvec\in \curlyC^1(\closOmT;\RR^3)$ with $\div\Svisc(\nablax\Uvec) \in \curlyC^0(\closOmT;\RR^3)$ and $C \in \curlyC^1(\closOmT)$ with $\nablax C \in \curlyC^1(\closOmT;\RR^3),\, \Deltax C \in \curlyC^1(\closOmT)$ are given satisfying
    \begin{align}
        \partial_t r + \div(r \Uvec) &= 0,\label{classical solution continuity}\\
        r\bigl( \partial_t \Uvec + \Uvec\cdot \nablax \Uvec\bigr) 
        &=
        \div\Svisc(\nablax \Uvec)
        -
        \nablax \Bigl( h(r) + q(r)\Bigr) 
        +
        \coup r \nablax \Bigl( C - r \Bigr),\label{classical solution momentum}\\
        \paracoup \partial_t C 
        -
        \kappa \Deltax C
        +
        \coup (C-r)
        &=
        0\label{classical solution parabolic}
    \end{align}
    pointwise in $\Omega_T$ and
    \begin{align}\label{classical solution BC}
        \Uvec_{\mid\partial\Omega} = 0,
        \quad 
        \nablax C\cdot \mathbf{n}_{\partial\Omega} = 0
    \end{align}
    pointwise in $[0,T]\times\partial\Omega$.
    The continuity equation $\eqref{classical solution continuity}$ implies that the renormalized continuity equation
    \begin{align}\label{classical solution renormalized continuity}
        \partial_t\bigl(b(r)\bigr)
        +
        \div\bigl( b(r) \Uvec\bigr)
        +
        \bigl( b^\prime(r) r - b(r) \bigr) \div\Uvec
        = 0
    \end{align}
    holds pointwise in $\Omega_T$ for any $b \in \curlyC^1(\RR)$.
    Applying the relative energy inequality $\eqref{proposition : relative energy inequality : inequality}$ from Proposition~\ref{proposition : relative energy inequality : inequality} for $(\rho,\uvec,c\mid r,\Uvec,C)$ yields for almost any $s \in (0,T)$ the inequality
    \begin{align}
        &\bigl[\Erel \bigr]_{\tau=0}^{\tau=s}
        +
        \int_0^s \int_\Omega 
        \lr{
        \Svisc(\nablax \uvec - \nablax \Uvec ) : (\nablax \uvec - \nablax \Uvec)
        +
        \paracoup|\partial_t c|^2 
        }
        \, \dd x \, \dd t\nonumber
        \\
        &\qquad \leq
        \mathcal{R}^1
        +
        \mathcal{R}^2
        +
        \mathcal{R}^3
        +
        \mathcal{R}^4,\label{WS uniqueness relative energy inequality I}
    \end{align}
    where $\mathcal{R}^1,\mathcal{R}^2,\mathcal{R}^3,\mathcal{R}^4$ are defined as in Proposition~\ref{proposition : relative energy inequality}.
    Our goal is to estimate the right-hand side in the relative energy inequality $\eqref{WS uniqueness relative energy inequality I}$ in terms of $\Erel$ in such a way that we can apply Gronwall's inequality later on.
    To do so, we use the equations $\eqref{classical solution continuity}$--$\eqref{classical solution renormalized continuity}$ that are satisfied by $(r,\Uvec,C)$.
    Let us treat each remainder separately.
    By using the momentum equation $\eqref{classical solution momentum}$, we obtain
    \begin{align}
        \mathcal{R}^1
        &=
        \int_0^s \int_\Omega 
        \rho \Bigl( r^{-1} \div\Svisc(\nablax\Uvec) - r^{-1}\nablax\Bigl( h(r) + q(r) \Bigr) + \coup \nablax\Bigl( C-r\Bigr)\Bigr)
        \cdot \Bigl(\Uvec-\uvec\Bigr)
        \, \dd x \, \dd t\nonumber
        \\
        &\quad 
        +
        \int_0^s \int_\Omega 
        \rho\Bigl(\uvec-\Uvec\Bigr) \cdot \nablax \Uvec \cdot \Bigl(\Uvec-\uvec\Bigr) \, \dd x \, \dd t\nonumber
        \\
        &\quad +\int_0^s \int_\Omega
        \biggr(
        \Svisc(\nablax\Uvec) :\Bigl(\nablax\Uvec - \nablax\uvec\Bigr) 
        +
        q(\rho) \Bigl( \div\uvec-\div\Uvec\Bigr)
        \biggr)
        \, \dd x \, \dd t\nonumber
        \\
        &\quad 
        +
        \int_0^s \int_\Omega
        \biggl(
        \Bigl( h(r) - h(\rho) \Bigr) \div\Uvec 
        +
        (r-\rho) \partial_t H^\prime(r)
        +
        \Bigl( r\Uvec - \rho \uvec\Bigr) \cdot \nablax  H^\prime(r) 
        \biggr)
        \, \dd x\, \dd t.\label{WS-Uniqueness REI prep:Term R^1 I}
    \end{align}
    Using integration by parts, we have
    \begin{equation*}
        \int_0^s \int_\Omega 
        \Svisc(\nablax \Uvec) : \nablax( \Uvec-\uvec )
        \, \dd x \, \dd t
        =
        -\int_0^s \int_\Omega 
        \div\Svisc(\nablax \Uvec) \cdot (\Uvec-\uvec)
        \, \dd x \, \dd t
    \end{equation*}
    and
    \begin{align*}
        &\int_0^s \int_\Omega 
        \frac{\rho}{r}\nablax \bigl( q(r) \bigr) \cdot \Bigl(\Uvec-\uvec\Bigr)
        \, \dd x \, \dd t
        \\
        &\quad 
        =
        \int_0^s \int_\Omega 
        \lr{
        \Bigl(\frac{\rho}{r} -1 \Bigr) \nablax q(r) \cdot \Bigl( \Uvec-\uvec\Bigr)
        +
        q(r)\Bigl( \div\uvec-\div\Uvec\Bigr) 
        }
        \,\dd x \, \dd t.
    \end{align*}
    By $\eqref{pressure pot:relations I}$ and $\eqref{pressure pot:relations II}$ we have
    \begin{align*}
        &\frac{\rho}{r}\nablax h (r)
        =
        \frac{\rho}{r}h^\prime(r) \nablax r
        =
        \rho H^{\prime\prime}(r)\nablax r
        =
        \rho \nablax H^\prime(r).
    \end{align*}
    By $\eqref{classical solution renormalized continuity}$, we infer
    \begin{align*}
        &\partial_t H^\prime(r)
        +
        \Uvec \cdot \nablax H^\prime(r)\\
        &=
        \partial_t H^\prime(r)
        +
        \div\bigl( H^\prime(r)\Uvec\bigr )
        +
        \bigl(H^{\prime\prime}(r) r - H^\prime(r) \bigr) \div\Uvec
        -
        H^{\prime\prime}(r) r \div\Uvec\\
        &=
        -H^{\prime\prime}(r) r \div\Uvec\\
        &=
        -h^\prime(r) \div\Uvec,
    \end{align*}
    which yields
    \begin{align*}
        &(r-\rho)\partial_t H^\prime(r)
        +
        \Bigl( r\Uvec - \rho \uvec \Bigr) \cdot \nablax  H^\prime(r)
        -
        \rho \nablax H^\prime(r) \cdot \Bigl( \Uvec - \uvec \Bigr) \\
        &=
        (r-\rho) \Bigl( \partial_t H^\prime(r) + \Uvec\cdot\nablax H^\prime(r) \Bigr)\\
        &=
        -(r-\rho)h^\prime(r)\div\Uvec.
    \end{align*}
    Using these identities in $\eqref{WS-Uniqueness REI prep:Term R^1 I}$ yields that 
    \begin{align*}
        \mathcal{R}^1
        =
        \mathcal{R}^{\mathrm{NSE}}
        +
        \mathcal{I}^1
        +
        \mathcal{I}^2
    \end{align*}
    with
    \begin{align*}
        \mathcal{R}^{\mathrm{NSE}}
        &:=
        \int_0^s \int_\Omega 
        \Bigl( \frac{\rho}{r}-1 \Bigr) \Bigl( \div\Svisc(\nablax\Uvec) - \nablax q(r)  \Bigr) \cdot \Bigl(\Uvec - \uvec\Bigr) 
        \, \dd x \, \dd t
        \\
        &\qquad 
        +
        \int_0^s \int_\Omega 
        \rho \Bigl(\uvec-\Uvec\Bigr) \cdot \nablax\Uvec \cdot \Bigl( \Uvec-\uvec\Bigr)
        \, \dd x \, \dd t\\
        &\qquad+
        \int_0^s \int_\Omega
        \Bigl( h(r) - h(\rho) - h^\prime(r) (r-\rho) \Bigr) \div\Uvec
        \, \dd x \, \dd t,
        \\
        \mathcal{I}^1
        &:=
        \int_0^s \int_\Omega 
        \Bigl( q(\rho) - q(r) \Bigr) \Bigl( \div\uvec - \div\Uvec\Bigr)
        \, \dd x \, \dd t,
        \\
        \mathcal{I}^2
        &:=
        \int_0^s \int_\Omega 
        \coup \rho \nablax C \cdot \Bigl( \Uvec- \uvec\Bigr) 
        -
        \coup \rho \nablax r \cdot \Bigl( \Uvec - \uvec \Bigr)
        \, \dd x \, \dd t.
    \end{align*}
    Notice that in view of an application of Gronwall's inequality, the term $\mathcal{R}^{\mathrm{NSE}}$ can be estimated appropriately by techniques that are available from the theory related to the compressible Navier--Stokes equations with constant viscosity coefficients (see e.g.~\cite{Feireisl2019, FeireislPetcuPravzak2019, Chaudhuri2019}). 
    These techniques work for any $\gamma \in (1,\infty)$.
    The term $\mathcal{I}^1$ can be estimated be by using the Lipschitz continuity of $q$ and Young's inequality as in \cite{Chaudhuri2019}, provided $\gamma \in [2,\infty)$.
    However, as we shall see in the sequel, we are able to estimate this term for the full range $\gamma \in (1,\infty)$ by exploiting the additional quadratic contributions in the relative energy functional that are not present for the compressible Navier--Stokes equations.
    Thus, we will not modify the terms $\mathcal{R}^{\mathrm{NSE}}$ and $\mathcal{I}^1$ further.
    We have now
    \begin{align}
        &[\Erel]_{\tau = 0}^{\tau = s}
        +
        \int_0^s \int_\Omega 
        \lr{
        \Svisc(\nablax \uvec - \nablax\Uvec) : (\nablax \uvec - \nablax\Uvec)
        +
        \paracoup|\partial_t c|^2 
        }
        \, \dd x \, \dd t\nonumber\\
        &\qquad \leq 
        \mathcal{R}^{\mathrm{NSE}} 
        +
        \mathcal{I}^1
        +
        \mathcal{I}^2
        +
        \mathcal{R}^2 + \mathcal{R}^3 + \mathcal{R}^4.\label{Modifications 0}
    \end{align}
    Let us now calculate the sum
    $
        \mathcal{I}^2
        +
        \mathcal{R}^2
        +
        \mathcal{R}^3
        +
        \mathcal{R}^4
    $
    more precisely by using the equations $\eqref{classical solution continuity}$, $\eqref{classical solution momentum}$ and $\eqref{classical solution parabolic}$.\\
    We have
    \begin{align}
        \mathcal{I}^2
        +
        \mathcal{R}^2\nonumber
        &=
        \int_0^s \int_\Omega 
        \lr{
        \coup\rho \nablax C \cdot \Uvec
        -
        \coup \rho \nablax r \cdot \Uvec
        -
        \coup \rho \nablax c \cdot \Uvec
        -
        \frac{\coup}{2} |\rho|^2 \div\Uvec
        }
        \, \dd x \, \dd t\nonumber\\
        &=
        \int_0^s \int_\Omega 
        \lr{
        \coup \rho \Bigl( \nablax C - \nablax c\Bigr) \cdot \Uvec
        -
        \coup \rho  \nablax r \cdot \Uvec
        -
        \frac{\coup}{2} |\rho|^2 \div \Uvec
        }
        \, \dd x \, \dd t.\label{Modifications I}
    \end{align}
    Then, we use the fact that $r$ satisfies the continuity equation $\eqref{classical solution continuity}$ to calculate $\mathcal{R}^3$ as
    \begin{align}
        \mathcal{R}^3
        &=
        -
        \int_0^s \int_\Omega 
        \lr{
        \coup
        (c-\rho) \div(r\Uvec)
        -
        \coup
        (C-r) \div(r\Uvec)
        }
        \, \dd x \, \dd t\nonumber
        \\
        &=
        \int_0^s \int_\Omega 
        \lr{
        \coup
        \rho \div(r\Uvec)
        -
        \coup
        c \div(r\Uvec)
        +
        \coup
        C \div(r\Uvec)
        -
        \coup
        r \div(r\Uvec) 
        }
        \, \dd x \, \dd t\nonumber
        \\
        &=
        \int_0^s \int_\Omega 
        \lr{
        \coup
        \rho \nablax r \cdot \Uvec 
        +
        \coup
        \rho r \div\Uvec
        +
        \coup
        r \nablax c \cdot \Uvec
        -
        \coup
        r \nablax C \cdot \Uvec
        +
        \coup
        r\nablax r \cdot \Uvec
        }
        \, \dd x \, \dd t\nonumber
        \\
        &=
        \int_0^s \int_\Omega 
        \lr{
        \coup
        \rho \nablax r \cdot \Uvec
        -
        \coup
        r\Bigl(\nablax C - \nablax c \Bigr) \cdot \Uvec
        +
        \coup
        \Bigl(\rho r - \frac{1}{2}|r|^2 \Bigr) \div\Uvec 
        }
        \, \dd x \, \dd t,\label{Modifications II}
    \end{align}
    where we have used integration by parts in the third and fourth step.
    By adding $\eqref{Modifications I}$ and $\eqref{Modifications II}$, we obtain
    \begin{align}
        &\mathcal{I}^2
        +
        \mathcal{R}^2
        +
        \mathcal{R}^3\nonumber
        \\
        &=
        \int_0^s \int_\Omega 
        \lr{
        \coup (\rho - r) \Bigl( \nablax C - \nablax c \Bigr) \cdot \Uvec
        -
        \frac{\coup}{2}
        \Bigl( |\rho|^2 - 2 \rho r + |r|^2 \Bigr)  \div\Uvec
        }
        \, \dd x \, \dd t\nonumber
        \\
        &=
        \int_0^s \int_\Omega 
        \lr{
        \coup (\rho-r)  \Bigl( \nablax C - \nablax c \Bigr) \cdot \Uvec
        -
        \frac{\coup}{2} |\rho - r|^2 \div\Uvec
        }
        \, \dd x \, \dd t.\label{Modifications III}
    \end{align}
    We further write by using integration by parts
    \begin{align}
        &\int_0^s \int_\Omega 
        \coup(\rho-r) \Bigl( \nablax C - \nablax c \Bigr) \cdot \Uvec
        \, \dd x \, \dd t\nonumber
        \\
        &\quad 
        =
        \int_0^s \int_\Omega 
        \biggl(
        \coup\Bigl( (\rho-r) - (c-C)\Bigr) \Bigl(\nablax C - \nablax c\Bigr) \cdot \Uvec
        +
        \coup (c-C)\Bigl(\nablax C - \nablax c\Bigr) \cdot \Uvec
        \biggr)
        \, \dd x \, \dd t\nonumber
        \\
        &\quad
        =
        \int_0^s \int_\Omega 
        \biggl(
        \coup \Bigl( (\rho-r) - (c-C)\Bigr) \Bigl( \nablax C - \nablax c \Bigr) \cdot \Uvec
        +
        \frac{\coup}{2} |c-C|^2 \div\Uvec
        \biggr)
        \, \dd x \, \dd t.\label{Modifications IV}
    \end{align}
    Now we calculate $\mathcal{R}^4$ with the help of the parabolic equation $\eqref{classical solution parabolic}$ that is satisfied by $C$ and also by $c$ in a strong sense.
    More specifically, we have
    \begin{align}
        \mathcal{R}^4
        =
        -
        \int_0^s \int_\Omega 
        \lr{
        \paracoup|\partial_t C|^2 
        -
        2 \paracoup \partial_t c \partial_t C
        }
        \, \dd x \, \dd t.\label{Modifications V}
    \end{align}
    Combining $\eqref{Modifications 0}$, $\eqref{Modifications III}$, $\eqref{Modifications IV}$ and $\eqref{Modifications V}$, we obtain
    \begin{align*}
        &\bigl[\Erel \bigr]_{\tau=0}^{\tau=s}
        +
        \int_0^s \int_\Omega 
        \lr{
        \Svisc(\nablax \uvec - \nablax \Uvec) : (\nablax \uvec - \nablax\Uvec)
        +
        \paracoup|\partial_t c-\partial_t C|^2
        }
        \, \dd x \, \dd t\nonumber
        \\
        &\quad
        \leq
        \mathcal{R}^{\mathrm{NSE}}
        +
        \mathcal{I}^1
        +
        \int_0^s \int_\Omega 
        \biggl(
        \coup \Bigl((\rho-r)-(c-C)\Bigr) \Bigl( \nablax C - \nablax c \Bigr) 
        \cdot \Uvec
        +
        \frac{\coup}{2} |c-C|^2 \div\Uvec
        \biggr)
        \, \dd x \, \dd t\nonumber
        \\
        &\qquad
        -
        \int_0^s \int_\Omega 
        \frac{\coup}{2} |\rho-r|^2 \div\Uvec
        \, \dd x \, \dd t.
    \end{align*}
    On our way proving Theorem~\ref{thm main result:W-S uniqueness}, we have shown the following partial result.
    \begin{proposition}\label{proposition : relative energy classical solution inserted}
        Let initial conditions $(\rho^0,(\rho\uvec)^0,c^0)$ be given satisfying $\eqref{defi FEW initial conditions}$ and let $(\rho,\uvec,c)$ be a finite energy weak solution existing on $[0,T]$ emanating from the initial conditions $(\rho^0,(\rho\uvec)^0,c^0)$.
        Let functions $r \in \curlyC^1(\closOmT)$ with $r>0$, $\Uvec\in \curlyC^1(\closOmT;\RR^3)$ with $\div \Svisc(\nablax\Uvec) \in \curlyC^0(\closOmT;\RR^3)$ and $C \in \curlyC^1(\closOmT)$ with $\nablax C \in \curlyC^1(\closOmT;\RR^3),\, \Deltax C \in \curlyC^1(\closOmT)$ be given satisfying $\eqref{classical solution continuity}$--$\eqref{classical solution BC}$.\\
        Then the inequality
        \begin{align}
            &[\Erel]_{\tau=0}^{\tau=s}
            +
            \int_0^s \int_\Omega 
            \lr{
            \Svisc(\nablax\uvec - \nablax\Uvec) : (\nablax \uvec - \nablax \Uvec)
            +
            \paracoup |\partial_t c - \partial_t C|^2 
            }
            \, \dd x \, \dd t\nonumber\\
            &\qquad \leq
            \mathcal{R}^{\mathrm{NSE}} 
            +
            \mathcal{J}^1
            +
            \mathcal{J}^2
            +
            \mathcal{J}^3\label{proposition : relative energy classical solution inserted : inequality}
        \end{align}
        holds for almost all $s \in (0,T)$, where
        \begin{align*}
            \mathcal{R}^{\mathrm{NSE}}
            &:=
            \int_0^s \int_\Omega 
            \Bigl( \frac{\rho}{r}-1 \Bigr) \Bigl( \div\Svisc(\nablax\Uvec) - \nablax q(r)  \Bigr) \cdot \Bigl(\Uvec - \uvec\Bigr) 
            \, \dd x \, \dd t\nonumber
            \\
            &\quad 
            +
            \int_0^s \int_\Omega 
            \biggl(
            \rho \Bigl(\uvec-\Uvec\Bigr) \cdot \nablax\Uvec \cdot \Bigl( \Uvec-\uvec\Bigr)
            +
            \Bigl( h(r) - h(\rho) - h^\prime(r) (r-\rho) \Bigr) \div\Uvec
            \biggr)
            \, \dd x \, \dd t,
            \\
            \mathcal{J}^1
            &:=
            \int_0^s \int_\Omega 
            \Bigl( q(\rho) - q(r) \Bigr) \Bigl( \div\uvec - \div\Uvec\Bigr) 
            \, \dd x \, \dd t,
            \\
            \mathcal{J}^2
            &:=
            \int_0^s \int_\Omega 
            \biggl(
            \coup \Bigl((\rho - r) - (c-C) \Bigr)\Bigl( \nablax C - \nablax c \Bigr) \cdot \Uvec
            +
            \frac{\coup}{2} |c-C|^2 \div\Uvec
            \biggr)
            \, \dd x \, \dd t,
            \\
            \mathcal{J}^3
            &:=
            -\int_0^s \int_\Omega 
            \frac{\coup}{2}
            |\rho-r|^2 
            \div\Uvec
            \, \dd x \, \dd t.
        \end{align*}
    \end{proposition}
    In order to estimate all the remainder on the right-hand side in inequality $\eqref{proposition : relative energy classical solution inserted : inequality}$, we need the following auxiliary lemma that guarantees a Poincaré inequality for the difference $c-C$ under a mild assumption on the initial data.
    Note that these assumption will hold in the proof of the weak-strong uniqueness property by hypothesis (cf. Theorem~\ref{thm main result:W-S uniqueness}).
    \begin{lemma}\label{lemma:WS uniqueness Poincare}
        Let the hypothesis and notation of Proposition~\ref{proposition : relative energy classical solution inserted} hold true and assume that
        \begin{align}\label{lemma:WS uniqueness Poincare IC}
            \int_\Omega 
            c^0
            \, \dd x
            =
            \int_\Omega 
            C(0)
            \, \dd x,
            \quad
            \int_\Omega 
            \rho^0
            \, \dd x
            =
            \int_\Omega
            r(0)
            \, \dd x.
        \end{align}
        Then we have for almost all $s \in (0,T)$ that
        \begin{align}\label{lemma:WS uniqueness Poincare inequality}
            ||c(s) - C(s)||_{L^2(\Omega)}
            \leq 
            C_P
            ||\nablax c(s) - \nablax C(s)||_{L^2(\Omega)},
        \end{align}
        where $C_P>0$ denotes the Poincaré constant on $\Omega$.
    \end{lemma}
    \begin{proof}
        For $s \in [0,T]$, we denote $\hat{c}(s):= c(s)-C(s)$ and $\hat{\rho}(s) := \rho(s)-r(s)$.
        Then $\hat{c}$ satisfies for any $s \in [0,T]$ and any test function $\varphi \in \curlyC^1([0,T]\times\overline{\Omega})$ the relation
        \begin{align*}
            \Bigl[\int_\Omega \paracoup\hat{c}(\tau)\varphi(\tau)\, \dd x  \Bigr]^{\tau=s}_{\tau=0}
            =
            \int_0^s
            \int_\Omega 
            \lr{
            \paracoup\hat{c} \partial_t \varphi
            -
            \kappa \nablax \hat{c} \cdot \nablax \varphi 
            -
            \coup(\hat{c}-\hat{\rho})\varphi
            }
            \, \dd x\, \dd t.
        \end{align*}
        By choosing $\varphi \equiv 1$ as a test function and using the first relation in $\eqref{lemma:WS uniqueness Poincare IC}$, we conclude for any $s \in [0,T]$ that
        \begin{align}
            \int_\Omega c(s) \, \dd x
            =
            -
            \frac{\coup}{\paracoup} 
            \int_0^s \int_\Omega 
            \hat{c}(\tau) 
            \, \dd x \, \dd \tau
            +
            \frac{\coup}{\paracoup} 
            \int_0^s \int_\Omega 
            \hat{\rho}(\tau)
            \, \dd x \, \dd t.\label{lemma:WS uniqueness Poincaré relation I}
        \end{align}
        By the continuity equation and the second relation in $\eqref{lemma:WS uniqueness Poincare IC}$, we have for any $s \in [0,T]$ that
        \begin{align*}
            \int_\Omega 
            \hat{\rho}(s)\, 
            \dd x
            =
            \int_\Omega 
            (\rho^0 - r(0))
            \, \dd x 
            =
            0.
        \end{align*}
        Thus, equation $\eqref{lemma:WS uniqueness Poincaré relation I}$ yields for any $s \in [0,T]$ that
        \begin{align*}
            \int_\Omega \hat{c}(s)
            \, \dd x 
            =
            -
            \frac{\coup}{\paracoup} 
            \int_0^s \int_\Omega 
            \hat{c}
            \, \dd x \, \dd t.
        \end{align*}
        With Gronwall's inequality, we conclude for any $s \in [0,T]$ that
        \begin{align*}
            \int_\Omega 
            \hat{c}(s) 
            \, \dd x
            = 0.
        \end{align*}
        With this relation, we deduce from Poincaré's inequality 
        \begin{align*}
            ||\hat{c}(s)||_{L^2(\Omega)}
            =
            ||\hat{c}(s) - \fint_\Omega \hat{c}(s)||_{L^2(\Omega)}
            \leq 
            C_P
            ||\nablax \hat{c}(s)||_{L^2(\Omega)}.
        \end{align*}
    \end{proof}

    We also need the following Lemma, in order to treat the term $\mathcal{R}^{\mathrm{NSE}}$ with techniques that are available from the theory on the compressible Navier--Stokes equations (see e.g.~\cite{FeireislJinNovotny2012,Feireisl2019,Chaudhuri2019,FeireislPetcuPravzak2019}).

    \begin{lemma}\label{lemma:auxiliary convexity estimates}
        Let $h \in C^0([0,\infty))\cap C^2((0,\infty))$ satisfy $\eqref{assump on p:monotonicity h}$ with $\gamma \in (1,\infty)$ and let $0<r_1<a<b<r_2 < \infty$ be given.
        Then, there exist two constants $k_h,K_h>0$ that only depend on $r_1,r_2,a,b$ and $h$, such that for any $r \in [a,b]$, we have the inequalities
        \begin{align}\label{lemma:auxiliary convexity estimates inequality I}
            H(\rho)
            -
            H(r)
            -
            H^\prime(r)(\rho-r)
            \geq 
            k_h
            \begin{cases}
                |\rho-r|^2 & \quad\rho \in [r_1,r_2]\\
                (1 + \rho^\gamma) & \quad\rho \in [0,r_1) \cup (r_2 ,\infty)
            \end{cases}
        \end{align}
        and
        \begin{align}\label{lemma:auxiliary convexity estimates inequality II}
            |h(\rho)-h(r)-h^\prime(r)(\rho-r)|
            \leq 
            K_h\Bigl( H(\rho)-H(r)-H^\prime(r)(\rho-r)\Bigr),
        \end{align}
        for any $\rho \in [0,\infty)$, where $H$ is defined as in $\eqref{pressure potential (def)}$.
    \end{lemma}
    \begin{proof}
        This follows from $\eqref{assump on p:monotonicity h}$ by using the relations in $\eqref{pressure pot:relations I}$, $\eqref{pressure pot:relations II}$ and Taylor's formula.
    \end{proof}

    We have now everything at hand to demonstrate the proof of our first main result Theorem~\ref{thm main result:W-S uniqueness}.
    \begin{proof}[Proof of Theorem~\ref{thm main result:W-S uniqueness}]
        We have that $\Erel_{\mid\tau=0} = 0$, and thus, by Proposition~\ref{proposition : relative energy classical solution inserted}, the inequality
        \begin{align}
            &\Erel(s)
            +
            \int_0^s \int_\Omega 
            \lr{
            \Svisc(\nablax \uvec - \nablax \Uvec) : (\nablax \uvec - \nablax \Uvec) 
            +
            \paracoup|\partial_t c - \partial_t C|^2 
            }
            \, \dd x \, \dd t\nonumber
            \\
            &\qquad\leq
            \mathcal{R}^{\mathrm{NSE}} 
            +
            \mathcal{J}^1
            +
            \mathcal{J}^2
            +
            \mathcal{J}^3,\label{Pf WS-uniqueness I}
        \end{align}
        which holds for almost all $s \in (0,T)$, where $\mathcal{R}^{\mathrm{NSE}},\mathcal{J}^1,\mathcal{J}^2$ and $\mathcal{J}^3$ are defined as in Proposition~\ref{proposition : relative energy classical solution inserted}.
        We estimate now the right-hand side of this inequality term-by-term. To shorten the notation we write 
        $\mathcal{E}:=\Erel$.
        The term $\mathcal{R}^{\mathrm{NSE}}$ can be estimated with the techniques that are available from the theory on the compressible Navier--Stokes equations (see e.g.~\cite{Feireisl2019, FeireislPetcuPravzak2019, Chaudhuri2019}).
        For the sake of completeness, we repeat these arguments here.
        We denote $\underline{r}:=\inf\limits_{(t,x)\in \closOmT} r(t,x)$, $\overline{r}:=\sup\limits_{(t,x) \in \closOmT} r(t,x)$, set $a:= \frac{1}{2}\underline{r}>0$, $b:= 2\overline{r}<\infty$ and take a cut-off function $\psi \in \curlyC^\infty_c((0,\infty))$ such that
        \begin{align*}
            0 \leq \psi \leq 1,\quad \psi = 1 \quad \text{on } [a,b].
        \end{align*}
        Then, we take $0<r_1<1<r_2$ such that $r_1<a$, $b < r_2$ and $\mathrm{supp}\,\psi \subseteq [r_1,r_2]$.
        For some integrable function $h \in L^1((0,T)\times\Omega)$, we decompose $h=h_{\mathrm{ess}}+h_{\mathrm{res}}$ with
        \begin{align*}
            h_{\mathrm{ess}}:=\psi(\rho) h,
            \quad 
            h_{\mathrm{res}}:= \bigl(1-\psi(\rho) \bigr) h.
        \end{align*}
        From inequality $\eqref{lemma:auxiliary convexity estimates inequality I}$ in Lemma~\ref{lemma:auxiliary convexity estimates} we conclude that there exists some positive constant $k_1>0$ such that
        \begin{align}\label{Pf WS-uniqueness II}
            \int_0^s\int_\Omega 
            \Bigl(
            \bigl[\uvec - \Uvec \bigr]^2_{\mathrm{ess}}
            +
            \bigl[\rho-r\bigr]_{\mathrm{ess}}^2
            +
            1_{\mathrm{res}}
            +
            \rho_{\mathrm{res}}
            +
            \rho^\gamma_{\mathrm{res}}
            \Bigr)
            \, \dd x \, \dd t
            \leq 
            k_1
            \int_0^s 
            \curlyE(t)
            \, \dd t.
        \end{align}
        For the first term in $\mathcal{R}^{\mathrm{NSE}}$ we first notice that by Rademacher's theorem, $q$ is almost everywhere differentiable in the classical sense and that the derivative is bounded by the Lipschitz constant $L_q \geq 0$ of $q$.
        This yields,
        \begin{align*}
            \|\nablax \bigl(q(r)\bigr)\|_{L^\infty(\OmegaT)}
            \leq 
            L_q \|\nablax r\|_{\curlyC(\closOmT)}
        \end{align*}
        and therefore
        \begin{align*}
            &\Biggl|
            \int_0^s \int_\Omega 
            \Bigl(\frac{\rho}{r}-1 \Bigr) \Bigl( \div\Svisc(\nablax \Uvec) - \nablax \bigl(q(r)\bigr)\Bigr) \cdot \Bigl( \Uvec - \uvec \Bigr)
            \, \dd x \, \dd t
            \Biggr|
            \\
            &\leq
            k_2
            \int_0^s \int_\Omega 
            \Bigl(
            |\rho-r|_{\mathrm{ess}} |\Uvec-\uvec|
            +
            |\rho-r|_{\mathrm{res}} |\Uvec-\uvec|
            \Bigr)
            \, \dd x \, \dd t
        \end{align*}
        with
        \begin{align*}
            k_2
            := 
            \rmin^{-1} \Bigl( \bigl\|\div\Svisc(\nablax\Uvec)\bigr\|_{\curlyC(\closOmT)} + L_q \|\nablax r\|_{\curlyC(\closOmT)} \Bigr).
        \end{align*}
        For $\delta \in (0,1)$, we obtain with Young's inequality and $\eqref{Pf WS-uniqueness II}$
        \begin{align*}
            k_2\int_0^s \int_\Omega 
            |\rho-r|_{\mathrm{ess}} |\Uvec - \uvec|
            \, \dd x \, \dd t
            &\leq 
            \frac{k_2^2}{4\delta}
            \int_0^s \int_\Omega 
            |\rho-r|^2_{\mathrm{ess}}
            \, \dd x \, \dd t
            +
            \delta 
            \int_0^s \int_\Omega 
            |\Uvec-\uvec|^2 
            \, \dd x \, \dd t
            \\
            &\leq 
            \frac{k_1k_2^2}{4\delta}
            \int_0^s 
            \curlyE(t)
            \, \dd t
            +
            \delta 
            \int_0^s 
            \int_\Omega 
            |\Uvec-\uvec|^2 
            \, \dd x,
        \end{align*}
        and
        \begin{align*}
            k_2&\int_0^s \int_\Omega 
            |\rho-r|_{\mathrm{res}} |\Uvec - \uvec |
            \, \dd x \, \dd t
            \\
            &\leq 
            k_2\int_0^s \int_\Omega 
            \Bigl(
            \rho_{\mathrm{res}}|\Uvec-\uvec|
            +
            r_{\mathrm{res}}|\Uvec-\uvec|
            \Bigr)
            \, \dd x \, \dd t
            \\
            &\leq 
            \int_0^s \int_\Omega 
            \Bigl(
            \rho_{\mathrm{res}}
            +
            k_2^2\,\rho_{\mathrm{res}}|\Uvec-\uvec|^2 
            +
            \frac{k_2^2\,\rmax^2}{4\delta} 1_{\mathrm{res}}
            +
            \delta |\Uvec-\uvec|^2 
            \Bigr)
            \, \dd x \, \dd t
            \\
            &\leq  
            \Bigl( k_1 + k_2^2 + \frac{k_1\,k_2^2\,\rmax^2}{4\delta} \Bigr)
            \int_0^s
            \curlyE(t)
            \, \dd t
            +
            \delta
            \int_0^s \int_\Omega 
            |\Uvec-\uvec|^2 
            \, \dd x \, \dd t.
        \end{align*}
        In total, we obtain
        \begin{align*}
            &\Biggl|
            \int_0^s \int_\Omega 
            \Bigl(\frac{\rho}{r}-1 \Bigr) \Bigl(\div\Svisc(\nablax \Uvec) - \nablax\bigl( q(r)\bigr) \Bigr) \cdot \Bigl( \Uvec - \uvec \Bigr)
            \, \dd x \, \dd t
            \Biggr|
            \\
            &\quad\leq 
            k_3(\delta) 
            \int_0^s 
            \curlyE(t)
            \, \dd t
            +
            2\delta 
            \int_0^s \int_\Omega 
            |\Uvec-\uvec|^2 
            \, \dd x \, \dd t
        \end{align*}
        with
        \begin{align*}
            k_3(\delta) := k_1 + k_2^2\Bigl(1 + \frac{k_1(1 + \rmax^2)}{4\delta}\Bigr).
        \end{align*}
        For the second term in $\mathcal{R}^{\mathrm{NSE}}$, we estimate with inequality $\eqref{lemma:auxiliary convexity estimates inequality II}$ from Lemma~\ref{lemma:auxiliary convexity estimates}
        \begin{align*}
            &\Biggl|
            \int_0^s \int_\Omega 
            \biggl(
            \rho\Bigl(\uvec-\Uvec\Bigr) \cdot \nablax \Uvec \cdot \Bigr( \Uvec-\uvec \Bigr)
            +
            \Bigl( h(r) - h(\rho) - h^\prime(r) (r-\rho) \Bigr) \div\Uvec
            \biggr)
            \, \dd x \, \dd t
            \Biggr|\\
            &\leq 
            k_4\int_0^s \mathcal{E}(t)\, \dd t.
        \end{align*}
        with
        \begin{align*}
            k_4 
            := 
            \|\nablax\Uvec\|_{\curlyC(\closOmT)} 
            +
            K_h\|\div\Uvec\|_{\curlyC(\closOmT)},
        \end{align*}
        where $K_h>0$ is a positive constant that only depends on $h$.
        Thus, we have
        \begin{align*}
            |\mathcal{R}^{\mathrm{NSE}}|
            \leq 
            k_5(\delta) 
            \int_0^s
            \mathcal{E}(t) 
            \, \dd t
            +
            2\delta 
            \int_0^s \int_\Omega 
            |\uvec-\Uvec|^2 
            \, \dd x \, \dd t,
        \end{align*}
        with $k_5(\delta) := k_3(\delta) + k_4$.
        Due to Poincaré's inequality and Korn's inequality, we obtain by choosing $\delta\in(0,1)$ small enough that
        \begin{align}\label{Pf WS-uniqueness III}
            |\mathcal{R}^{\mathrm{NSE}}| 
            \leq 
            k_6 
            \int_0^s 
            \mathcal{E}(t)
            \, \dd t
            +
            \frac{1}{4} 
            \int_0^s \int_\Omega 
            \Svisc(\nablax\uvec-\nablax\Uvec):(\nablax\uvec-\nablax \Uvec)
            \, \dd x \, \dd t,
        \end{align}
        for some positive constant $k_6>0$.\\
        For $\mathcal{J}^1$, we first observe that due to the assumption $\rho^0=r(0)$ and $c^0=C(0)$ the conditions $\eqref{lemma:WS uniqueness Poincare IC}$ are fulfilled.
        Thus, we are allowed to apply inequality $\eqref{lemma:WS uniqueness Poincare inequality}$ from Lemma~$\ref{lemma:WS uniqueness Poincare}$ and obtain with the help of Young's inequality
        \begin{align}
            \int_0^s \int_\Omega 
            |\rho-r|^2 
            \, \dd x \, \dd t\nonumber
            &\leq 
            \int_0^s \int_\Omega 
            \lr{
            2| (\rho-r) - (c-C)|^2 
            +
            2|c-C|^2 
            }
            \, \dd x \, \dd t\nonumber
            \\
            &\leq 
            \int_0^s \int_\Omega 
            \lr{
            2 |(\rho-r) - (c-C)|^2 
            +
            2C_P^2 |\nablax c - \nablax C|^2 
            }
            \, \dd x \, \dd t\nonumber
            \\
            &\leq 
            4\Bigl( \frac{1}{\coup} + \frac{ C_P^2}{\kappa}\Bigr)
            \int_0^s 
            \mathcal{E}(t)
            \, \dd t.\label{Pf WS-uniqueness IV}
        \end{align}
        With the Lipschitz continuity of $q$ and another application of Young's inequality, we obtain
        \begin{align}
            |\mathcal{J}^1|
            &\leq 
            \frac{1}{\lambda} 
            \int_0^s \int_\Omega 
            | q(\rho) - q(r)|^2 
            \, \dd x \, \dd t
            +
            \frac{\lambda}{4}
            \int_0^s \int_\Omega 
            |\div \uvec - \div \Uvec|^2 
            \, \dd x \, \dd t\nonumber
            \\
            &\leq 
            \frac{L_q^2}{\lambda}
            \int_0^s \int_\Omega |\rho - r|^2 
            \, \dd x \, \dd t
            +
            \frac{1}{4}
            \int_0^s \int_\Omega 
            \Svisc(\nablax\uvec-\nablax\Uvec):(\nablax\uvec-\nablax\Uvec)
            \, \dd x \, \dd t\nonumber
            \\
            &\leq 
            k_7
            \int_0^s 
            \mathcal{E}(t)
            \, \dd t
            +
            \frac{1}{4} 
            \int_0^s \int_\Omega 
            \Svisc(\nablax\uvec-\nablax\Uvec):(\nablax\uvec - \nablax\Uvec)
            \, \dd x \, \dd t\label{Pf WS-uniqueness V}
        \end{align}
        with
        \begin{align*}
            k_7 
            :=
            \frac{4L_q^2}{\lambda}\Bigl(\frac{1}{\coup} + \frac{C_P^2}{\kappa} \Bigr).
        \end{align*}
        For $\mathcal{J}^2$, we obtain with Young's inequality and inequality $\eqref{lemma:WS uniqueness Poincare inequality}$ from Lemma~\ref{lemma:WS uniqueness Poincare}
        \begin{align*}
            |\mathcal{J}^2|
            &\leq 
            \|\Uvec\|_{\curlyC(\closOmT)}
            \int_0^s \int_\Omega 
            \frac{\coup}{2} | (\rho-r) - (c-C)|^2 
            +
            \frac{\coup}{2} |\nablax c - \nablax C|^2 
            \, \dd x \, \dd t
            \\
            &\qquad 
            +
            \|\div\Uvec\|_{\curlyC(\closOmT)}
            \int_0^s \int_\Omega 
            \frac{\coup C_P^2}{2}|\nablax c - \nablax C|^2 
            \, \dd x \, \dd t
        \end{align*}
        and for $\mathcal{J}^3$ we obtain by using $\eqref{Pf WS-uniqueness IV}$
        \begin{align*}
            |\mathcal{J}_3|
            \leq 
            \|\div\Uvec\|_{\curlyC(\closOmT)}
            \int_0^s \int_\Omega
            \frac{\coup}{2} |\rho-r|^2
            \, \dd x \, \dd t
            \leq 
            2\|\div\Uvec\|_{\curlyC(\closOmT)}\Bigl(1+\frac{\coup C_P^2}{\kappa}\Bigr) 
            \int_0^s 
            \mathcal{E}(t)
            \, \dd t.
        \end{align*}
        This yields
        \begin{align}\label{Pf WS-uniqueness VI}
            |\mathcal{J}^2| + |\mathcal{J}^3| 
            \leq 
            k_8 
            \int_0^s
            \mathcal{E}(t) 
            \, \dd t
        \end{align}
        with
        \begin{align*}
            k_8 
            := 
            \|\Uvec\|_{\curlyC(\closOmT)} \Bigl(1 + \frac{\coup}{\kappa}\Bigr) 
            +
            \|\div\Uvec\|_{\curlyC(\closOmT)} \Bigl( 2 +  \frac{3\coup C_P^2}{\kappa}\Bigr).
        \end{align*}
        Combining the inequalities $\eqref{Pf WS-uniqueness I}$, $\eqref{Pf WS-uniqueness III}$, $\eqref{Pf WS-uniqueness V}$ and $\eqref{Pf WS-uniqueness VI}$ implies that the inequality
        \begin{align*}
             \mathcal{E}(s) 
             \leq 
             k_{9}
             \int_0^s \mathcal{E}(t) 
             \, \dd t
        \end{align*}
        holds for almost all $s \in (0,T)$, where $k_{9}:= k_6 + k_7 + k_8$.
        Applying Gronwall's inequality yields that $\mathcal{E}(s) = 0$ for almost all $s \in (0,T)$ and the result follows.
    \end{proof}
    
    \section{The Relaxation Limit}\label{Sec:Singular Limit}
    This section is devoted to the proof of Theorem~\ref{thm main result:singular limit}.
    More precisely, for $\coup>0$ we consider the initial-boundary value problem $\eqref{relaxed NSK field eq continuity}$--$\eqref{relaxed NSK IC}$ with $\paracoup=\paracoup(\coup)$ being a function of $\coup$ with $\paracoup(\coup)\to 0$ as $\coup \to \infty$ and prove that a sequence of finite energy weak solutions to $\eqref{relaxed NSK field eq continuity}$--$\eqref{relaxed NSK IC}$ converges in certain norms to a smooth solution of the NSK equations $\eqref{NSK field eq continuity}$--$\eqref{NSK IC}$, provided such a solution exists and the initial data is not too ill-prepared.
    In the course of the proof, we will identify a convergence rate for the corresponding norms as a byproduct.
    As a tool, we will use again the relative energy inequality derived in Section~\ref{Sec:Relative Energy Inequality}, however, this time, we will test this inequality with a classical solution to the target system, that is, the compressible NSK equations.
    More specifically, for $\coup,\paracoup>0$, we fix initial data $(\rho_\coup^0,(\rho\uvec)_\coup^0,c_\coup^0)$ satisfying $\eqref{defi FEW initial conditions}$ and a finite energy weak solution $(\rho_\coup,\uvec_\coup,c_\coup)$ of $\eqref{relaxed NSK field eq continuity}$--$\eqref{relaxed NSK IC}$ existing on $[0,T]$ emanating from the initial data $(\rho_\coup^0,(\rho\uvec)_\coup^0,c_\coup^0)$.
    As we will assume later in the proof of Theorem~\ref{thm main result:singular limit} that $\paracoup=\paracoup(\coup)$ is a function of $\coup>0$ with $\paracoup(\coup)\to 0$ for $\coup \to \infty$, we omit the dependency on $\paracoup$ in the notation.
    We assume that functions $r\in \curlyC^1(\closOmT)$ with $r>0, \, \nablax r \in \curlyC^1(\closOmT;\RR^3),\, \Deltax r \in \curlyC^1(\closOmT)$ and $\Uvec\in \curlyC^1(\closOmT;\RR^3)$ with $\div\Svisc(\nablax\Uvec), \, \nablax \div\Uvec \in \curlyC^0(\closOmT;\RR^3)$ are given satisfying
    \begin{align}
        \partial_t r 
        +
        \div(r \Uvec)
        &= 
        0,\label{classical solution NSK continuity}
        \\
        r\bigl( \partial_t\Uvec + \Uvec\cdot \nablax\Uvec\bigr)
        &=
        \div\Svisc(\nablax\Uvec)
        -
        \nablax p(r)
        +
        \kappa r \Deltax\nablax r\label{classical solution NSK momentum}
    \end{align}
    pointwise in $\Omega_T$ and
    \begin{align}\label{classical solution NSK BC}
        \Uvec_{\mid\partial\Omega} 
        =
        0,
        \quad
        \nablax r \cdot \mathbf{n}_{\mid\partial\Omega}
        =
        0
    \end{align}
    pointwise in $[0,T]\times\partial\Omega$. The continuity equation $\eqref{classical solution NSK continuity}$ implies that the renormalized continuity equation $\eqref{classical solution renormalized continuity}$ also holds pointwise.
    Since $r$ satisfies a zero Neumann boundary condition, we can apply Proposition~\ref{proposition : relative energy inequality} with $C=r$ and infer that the relative energy inequality
    \begin{align*}
        &\bigl[ \mathcal{E}(\rho_\coup,\uvec_\coup,c_\coup\mid r, \Uvec) \bigr]_{\tau=0}^{\tau=s}
        +
        \int_0^s \int_\Omega
        \lr{
        \Svisc(\nablax \auvec-\nablax\Uvec):(\nablax\auvec - \nablax \Uvec)
        +
        \paracoup|\partial_t \ac|
        }
        \, \dd x \, \dd t\nonumber
        \\
        &\qquad \leq 
        \mathcal{R}^1_\coup
        +
        \mathcal{R}^2_\coup
        +
        \mathcal{R}^3_\coup
        +
        \mathcal{R}^4_\coup
    \end{align*}
    holds for almost all $s \in (0,T)$, where  
    \begin{align}
        \mathcal{E}(\arho,\auvec,\ac\mid r,\Uvec) 
        :&=
        \mathcal{E}(\arho,\ac,\auvec \mid r,\Uvec,C=r)\nonumber
        \\
        &=
        \int_\Omega 
        \lr{
        \frac{1}{2} \arho |\auvec-\Uvec|^2 
        +
        H(\arho)-H(r)-H^\prime(r)(\arho-r)
        }
        \, \dd x\nonumber
        \\
        &\qquad+
        \int_\Omega
        \lr{
        \frac{\coup}{2}|\arho-\ac|^2 
        +
        \frac{\kappa}{2}|\nablax \ac - \nablax r|^2
        }
        \, \dd x,\label{defi Erel for relax limit}
    \end{align}
    and the remainders $\mathcal{R}^1_\coup,\mathcal{R}^2_\coup,\mathcal{R}^3_\coup,\mathcal{R}^4_\coup$ are given by
    \begin{align*}
        \mathcal{R}^1_\coup 
        &:= 
        \tilde{\mathcal{R}}^{\mathrm{NSE}}_{\coup} + \tilde{\mathcal{R}}^1_\coup,
        \\
        \tilde{\mathcal{R}}^{\mathrm{NSE}}_{\coup} 
        &:=
        \int_0^s \int_\Omega 
        \Bigl( \frac{\arho}{r} - 1\Bigr) \Bigl(\div\Svisc(\nablax\Uvec) - \nablax \bigl( q(r)\bigr) \Bigr) \cdot \Bigl( \Uvec-\auvec \Bigr)
        \, \dd x \, \dd t,
        \\
        &\qquad+
        \int_0^s \int_\Omega 
        \arho\Bigl(\auvec-\Uvec\Bigr) \cdot \nablax \Uvec \cdot \Bigl( \Uvec-\auvec\Bigr)
        \, \dd x \, \dd t
        \\
        &\qquad+
        \int_0^s \int_\Omega
        \Bigl( h(r) - h(\arho) - h^\prime(r)(r-\arho) \Bigr) \div\Uvec
        \, \dd x\, \dd t
        \\
        &\qquad+
        \int_0^s \int_\Omega 
        \Bigl( q(\arho) - q(r) \Bigr) \Bigl( \div\auvec - \div\Uvec \Bigr)
        \, \dd x\, \dd t,
        \\
        \tilde{\mathcal{R}}^1_\coup 
        &:=
        \int_0^s \int_\Omega 
        \kappa \arho \nablax \Deltax r \cdot \Bigl( \Uvec-\auvec \Bigr)
        \, \dd x\, \dd t,
        \\
        \mathcal{R}^2_\coup
        &:=
        -
        \int_0^s \int_\Omega 
        \lr{
        \frac{\coup}{2} |\arho|^2 \div\Uvec
        +
        \coup \arho\nablax \ac \cdot \Uvec
        }
        \, \dd x\, \dd t,
        \\
        \mathcal{R}^3_\coup
        &:=
        \int_0^s \int_\Omega 
        \coup (\ac - \arho) \partial_t r
        \, \dd x\, \dd t,
        \\
        \mathcal{R}^4_\coup
        &:=
        \int_0^s \int_\Omega 
        \lr{
        \kappa \Deltax \ac \partial_t  r
        -
        \kappa \Deltax r \partial_t r
        -
        \coup (\ac - \arho) \partial_t r
        +
        \kappa \Deltax r \partial_t \ac
        }
        \, \dd x\, \dd t.
    \end{align*}
    By using the continuity equation $\eqref{classical solution NSK continuity}$, we rewrite
    \begin{align*}
        \tilde{\mathcal{R}}^1_\coup
        +
        \mathcal{R}^2_\coup
        +
        \mathcal{R}^3_\coup
        +
        \mathcal{R}^4_\coup
        =
        \mathcal{I}^1_\coup
        +
        \mathcal{I}^2_\coup,
    \end{align*}
    with
    \begin{align*}
        \mathcal{I}^1_\coup
        &:=
        \int_0^s \int_\Omega
        \lr{
        \kappa \arho \nablax \Deltax r \cdot \Uvec
        -
        \frac{\coup}{2}|\arho|^2 \div\Uvec
        -
        \coup \arho \nablax \ac \cdot \Uvec
        }
        \, \dd x \, \dd t
        \\
        &\quad+
        \int_0^s \int_\Omega
        \lr{
        \kappa \Deltax r \nablax r \cdot \Uvec
        +
        \kappa r \Deltax r \div\Uvec
        -
        \kappa \Deltax \ac \nablax r \cdot \Uvec
        -
        \kappa r \Deltax \ac \div\Uvec
        }
        \, \dd x \, \dd t,
        \\
        \mathcal{I}_\coup^2
        &:=
        \int_0^s \int_\Omega 
        \lr{
        \kappa \Deltax r \partial_t \ac 
        -
        \kappa \arho \nablax \Deltax r \cdot \auvec
        }
        \, \dd x \, \dd t.
    \end{align*}
    We have to manipulate the term $\mathcal{I}_\coup^1$ further, so that we can estimate this term appropriately later on.
    In fact, we claim that the following identity holds
    \begin{align}
        &\mathcal{I}_\coup^1
        +
        \int_0^s \int_\Omega 
        \lr{
        \paracoup\partial_t\ac\div(\ac \Uvec)
        +
        \frac{\coup}{2}|\arho-\ac|^2 \div
        \Uvec
        }
        \, \dd x \, \dd t\nonumber
        \\
        &=
        \int_0^s \int_\Omega
        \lr{
        \kappa \Bigl((r-\ac) \Deltax (r-\ac) + \frac{|\nablax(r-\ac)|^2}{2}\Bigr)\div\Uvec
        }
        \, \dd x \, \dd t\nonumber
        \\
        &\quad-
        \int_0^s \int_\Omega
        \biggl(
        \kappa \Bigl( \nablax(r-\ac)\otimes \nablax (r-\ac) \Bigr):\nablax \Uvec
        +
        \kappa (\ac-\arho) \nablax\Deltax r \cdot \Uvec
        \biggr)
        \, \dd x\, \dd t.\label{singular limit decomposition I_coup^1 I}
    \end{align}
    To see this, we compute the right-hand side of this identity term-by-term. By using integration by parts, we infer that
    \begin{align*}
        \mathcal{T}^1_\coup
        &:=
        \int_0^s \int_\Omega 
        \kappa (r-\ac )\Deltax (r-\ac ) \div \Uvec
        \, \dd x \, \dd t\\
        &=
        \int_0^s\int_\Omega 
        \lr{
        \kappa r  \Deltax r \div\Uvec
        +
        \kappa \ac \Deltax \ac \div\Uvec
        -
        \kappa r  \Deltax \ac \div\Uvec
        -
        \kappa \ac \Deltax r \div \Uvec
        }
        \, \dd x \, \dd tr,
        \\
        \mathcal{T}^2_\coup
        &:=
        \int_0^s\int_\Omega 
        \lr{
        \frac{\kappa}{2}  |\nablax(r-\ac)|^2 \div\Uvec
        -
        \kappa \Bigl[ \nablax(r-\ac)\otimes\nablax(r-\ac) \Bigr] : \nablax \Uvec
        }
        \, \dd x \, \dd t
        \\
        &=\int_0^s \int_\Omega 
        \lr{
        \kappa \Deltax r \nablax r \cdot \Uvec
        +
        \kappa \Deltax \ac \nablax \ac \cdot \Uvec
        -
        \kappa \Deltax r \nablax \ac \cdot \Uvec
        -
        \kappa \Deltax \ac \nablax r \cdot \Uvec
        }
        \, \dd x \, \dd t,
        \\
        \mathcal{T}^3_\coup
        &:=
        -\int_0^s \int_\Omega 
        \kappa (\ac - \arho) \nablax \Deltax r \cdot \Uvec
        \, \dd x \, \dd t\\
        &=
        -\int_0^s \int_\Omega 
        \lr{
        \kappa \ac \nablax \Deltax r \cdot \Uvec
        -
        \kappa \arho \nablax \Deltax r \cdot \Uvec
        }
        \, \dd x \, \dd t,
    \end{align*}
    Thus, several terms cancel out, if we subtract $\mathcal{T}^1_\coup+\mathcal{T}_\coup^2+\mathcal{T}_\coup^3$ from $\mathcal{I}_\coup^1$.
    We obtain
    \begin{align}
        &\mathcal{I}_\coup^1 
        -
        \mathcal{T}^1_\coup
        -
        \mathcal{T}^2_\coup
        -
        \mathcal{T}^3_\coup\nonumber 
        \\
        &=
        \int_0^s \int_\Omega 
        \lr{
        \kappa \ac \Deltax r \div\Uvec
        -
        \kappa \Deltax \ac \div(\ac\Uvec)
        +
        \kappa \Deltax r \nablax \ac \cdot \Uvec
        +
        \kappa \ac \nablax \Deltax r \cdot \Uvec
        }
        \, \dd x \, \dd t\nonumber
        \\
        &\quad
        -
        \int_0^s \int_\Omega
        \lr{
        \frac{\coup}{2}|\arho|^2 \div\Uvec
        +
        \coup \arho\nablax \ac \cdot\Uvec
        }
        \, \dd x\, \dd t\nonumber
        \\
        &=
        -
        \int_0^s \int_\Omega 
        \lr{
        \kappa \Deltax \ac \div(\ac \Uvec)
        +
        \frac{\coup}{2}|\arho|^2 \div\Uvec
        +
        \coup\arho \nablax\ac \cdot \Uvec
        }
        \, \dd x\, \dd t,\label{singular limit decomposition I_coup^1 II}
    \end{align}
    where we have used integration by parts to obtain the second equality.
    By using the parabolic equation 
    \begin{align*}
        \paracoup\partial_t \ac 
        -
        \kappa \Deltax \ac 
        =
        \coup(\arho -\ac)
    \end{align*}
    that holds pointwise almost everywhere in $(0,T)\times\Omega$, we compute further
    \begin{align}\label{Repair I}
        -\int_0^s \int_\Omega 
        \kappa \Deltax \ac \div(\ac \Uvec)
        \, \dd x \, \dd t
        =
        -
        \int_0^s \int_\Omega 
        \Bigl(
        \paracoup\partial_t \ac \div(\ac \Uvec) 
        -
        \coup(\arho-\ac)\div(\ac \Uvec)
        \Bigr)
        \, \dd x \, \dd t
    \end{align}
    and
    \begin{align}
        &\int_0^s \int_\Omega 
        \coup(\arho -\ac) \div(\ac \Uvec)
        \, \dd x \, \dd t\nonumber
        \\
        &=
        \int_0^s \int_\Omega 
        \lr{
        \coup\arho \nablax \ac \cdot \Uvec
        +
        \coup\arho \ac \div\Uvec
        -
        \coup \ac \nablax \ac \cdot \Uvec
        -
        \coup |\ac |^2 \div\Uvec
        }
        \, \dd x\, \dd t\nonumber
        \\
        &=
        \int_0^s \int_\Omega 
        \lr{
        \coup\arho\nablax \ac \cdot \Uvec
        +
        \coup\arho \ac \div\Uvec
        -
        \frac{\coup}{2}|\ac|^2 \div
        \Uvec
        }
        \, \dd x \, \dd t\nonumber
        \\
        &=
        \int_0^s \int_\Omega 
        \lr{
        \coup\arho\nablax\ac \cdot \Uvec
        -
        \frac{\coup}{2}|\arho-\ac|^2 \div\Uvec
        +
        \frac{\coup}{2}|\arho|^2 \div\Uvec
        }
        \, \dd x\, \dd t.\label{Repair II}
    \end{align}
    Substituting $\eqref{Repair I}$ and $\eqref{Repair II}$ into $\eqref{singular limit decomposition I_coup^1 II}$ yields
    \begin{align*}
        \mathcal{I}_\coup^1 
        -
        \mathcal{T}^1_\coup
        -
        \mathcal{T}^2_\coup
        -
        \mathcal{T}^3_\coup
        =
        -
        \int_0^s \int_\Omega 
        \lr{
        \frac{\coup}{2}|\arho-\ac|^2 \div\Uvec
        +
        \paracoup\partial_t \ac \div(\ac\Uvec)
        }
        \, \dd x\, \dd t,
    \end{align*}
    which is exactly $\eqref{singular limit decomposition I_coup^1 I}$.
    On our way proving Theorem~\ref{thm main result:singular limit}, we have shown the following partial result.
    \begin{proposition}\label{proposition : REI singular limit}
        For $\coup,\paracoup>0$ let initial conditions $(\rho_\coup^0,(\rho\uvec)_\coup^0,c_\coup^0)$ be given satisfying $\eqref{defi FEW initial conditions}$.
        Let $(\rho_\coup,\uvec_\coup,c_\coup)$ denote a finite energy weak solution of $\eqref{relaxed NSK field eq continuity}$--$\eqref{relaxed NSK IC}$ existing on $[0,T]$ emanating from the initial conditions $(\rho^0_\coup,(\rho\uvec)^0_\coup,c^0_\coup)$.
        Let functions $r \in \curlyC^1(\closOmT)$ with $r>0, \, \nablax r \in \curlyC^1(\closOmT;\RR^3),\, \Deltax r \in \curlyC^1(\closOmT)$ and $\Uvec\in \curlyC^1(\closOmT;\RR^3)$ with $\div\Svisc(\nablax\Uvec),\, \nablax\div\Uvec \in \curlyC^0(\closOmT;\RR^3)$ be given satisfying $\eqref{classical solution NSK continuity}$--$\eqref{classical solution NSK BC}$.\\
        Then the relative energy inequality
        \begin{align}
            &\bigl[\curlyE(\arho,\auvec,\ac\mid r,\Uvec)\bigr]_{\tau=0}^{\tau=s}
            +
            \int_0^s \int_\Omega 
            \lr{
            \Svisc(\nablax\auvec-\nablax\Uvec):(\nablax\auvec-\nablax\Uvec)
            +
            \paracoup|\partial_t\ac|^2 
            }
            \, \dd x \, \dd t\nonumber
            \\
            &\qquad 
            \leq 
            \mathcal{R}^{\mathrm{NSE}}_\coup
            +
            \sum\limits_{i=0}^{7} 
            \mathcal{J}^i_\coup\label{proposition : REI singular limit inequality}
        \end{align}
        holds for almost all $s \in (0,T)$, where $\curlyE(\arho,\auvec,\ac\mid r,\Uvec)$ is defined as in $\eqref{defi Erel for relax limit}$ and
        \begin{align*}
            \mathcal{R}^{\mathrm{NSE}}_\coup
            &:=
            \int_0^s \int_\Omega 
            \Bigl( \frac{\arho}{r} - 1\Bigr) \Bigl( \div\Svisc(\nablax\Uvec) - \nablax\bigl(q(r)\bigr)\Bigr) \cdot \Bigl( \Uvec-\auvec \Bigr)
            \, \dd x \, \dd t\nonumber
            \\
            &\qquad+
            \int_0^s \int_\Omega 
            \arho\Bigl(\auvec-\Uvec\Bigr) \cdot \nablax \Uvec \cdot \Bigl( \Uvec-\auvec\Bigr)
            \, \dd x \, \dd t
            \\
            &\qquad+
            \int_0^s \int_\Omega
            \Bigl( h(r) - h(\arho) - h^\prime(r) (r-\arho)\Bigr) \div\Uvec
            \, \dd x\, \dd t,
            \\
            \mathcal{J}^1_\coup
            &:=
            \int_0^s \int_\Omega 
            \Bigl( q(\arho) - q(r) \Bigr) \Bigl( \div\auvec - \div\Uvec\Bigr)
            \, \dd x\, \dd t,
            \\
            \mathcal{J}^2_\coup
            &:=
            \int_0^s\int_\Omega 
            \kappa \Bigl((r-\ac) \Deltax (r-\ac) + \frac{|\nablax(r-\ac)|^2}{2}\Bigr) \div\Uvec
            \, \dd x\, \dd t
            \\
            \mathcal{J}^3_\coup
            &:=
            -
            \int_0^s \int_\Omega 
            \kappa \Bigl( \nablax(r-\ac)\otimes \nablax (r-\ac) \Bigr):\nablax \Uvec
            \, \dd x\, \dd t
            \\
            \mathcal{J}^4_\coup
            &:=
            -
            \int_0^s \int_\Omega 
            \kappa (\ac-\arho) \Uvec \cdot \nablax\Deltax r
            \, \dd x \, \dd t
            \\
            \mathcal{J}^5_\coup
            &:=
            -
            \int_0^s \int_\Omega
            \frac{\coup}{2}|\arho-\ac|^2 \div\Uvec
            \, \dd x\, \dd t
            \\
            \mathcal{J}^6_\coup
            &:=
            -
            \int_0^s \int_\Omega
            \paracoup \partial_t\ac \div(\ac\Uvec)
            \, \dd x\, \dd t
            \\
            \mathcal{J}^7_\coup 
            &:=
            \int_0^s \int_\Omega
            \lr{
            \kappa \Deltax r
            \partial_t\ac 
            -
            \kappa\arho\auvec \cdot \nablax \Deltax r
            }
            \, \dd x \, \dd t.
        \end{align*}
    \end{proposition}
    In order to estimate the remainders appearing on the right-hand side in the relative energy inequality $\eqref{proposition : REI singular limit inequality}$ appropriately, we need the following auxiliary lemma that follows from Poincaré's inequality.    
    \begin{lemma}\label{lemma : generalized Poincaré - singular limit}
        Let the hypothesis of Proposition~\ref{proposition : REI singular limit} hold true.
        Then we have for any $\coup,\paracoup>0$ and almost all $s \in (0,T)$ the estimate
        \begin{align}\label{lemma : generalized Poincaré - singular limit : inequality}
            \|r(s)-\ac(s)\|_{L^2(\Omega)}
            \leq 
            C_p \|\nablax r(s) - \nablax \ac(s)\|_{L^2(\Omega)}
            +
            \frac{\paracoup}{\coup}
            \|\partial_t \ac (s) \|_{L^2(\Omega)}
            +
            |e_\coup| \sqrt{|\Omega|},
        \end{align}
        where $C_P>0$ denotes the Poincaré constant on $\Omega$ and
        \begin{align*}
            e_\coup 
            := 
            \fint_\Omega (\arho^0 - r(0)) \, \dd x.
        \end{align*}
    \end{lemma}
    \begin{proof}
        We fix $\coup>0$ and use that $\ac$ satisfies the linear parabolic equation
        \begin{align*}
            \paracoup \partial_t \ac 
            -
            \kappa \Deltax \ac
            = 
            \coup (\arho - \ac)
        \end{align*}
        almost everywhere in $\Omega_T$ with zero Neumann boundary conditions in the sense of traces.
        This yields for almost any $s \in (0,T)$ by taking the average over $\Omega$ and multiplying by $\coup^{-1}$
        \begin{align}\label{lemma : generalized Poincaré - singular limit : I}
            \frac{\paracoup}{\coup}
            \fint_\Omega 
            \partial_t \ac (s)
            \, \dd x 
            =
            \fint_\Omega
            (\arho-\ac)(s)
            \, \dd x.
        \end{align}
        We conclude with the continuity equation for $r$ and $\arho$ that
        \begin{align*}
            \int_\Omega 
            \arho(s)
            \, \dd x 
            =
            \int_\Omega
            \arho^0
            \, \dd x
            ,
            \quad 
            \int_\Omega
            r(0)
            \, \dd x
            =
            \int_\Omega 
            r(s) 
            \, \dd x.
        \end{align*}
        Thus,
        \begin{align*}
            \fint_\Omega (\arho-r)(s)
            \, \dd x 
            = 
            \fint_\Omega 
            (\arho^0 - r(0))
            \, \dd x
            = e_\coup.
        \end{align*}
        With this relation, we conclude from $\eqref{lemma : generalized Poincaré - singular limit : I}$ that
        \begin{align}\label{lemma : generalized Poincaré - singular limit : II}
            \frac{\paracoup}{\coup}
            \fint_\Omega 
            \partial_t \ac(s)
            \, \dd x
            =
            \fint_\Omega 
            (r-\ac)(s)
            \, \dd x
            +
            e_\coup.
        \end{align}
        Now we make use of the Poincaré inequality
        \begin{align*}
            \Bigl\|\phi - \fint_\Omega \phi\Bigl\|_{L^2(\Omega)}
            \leq 
            C_P||\nablax \phi||_{L^2(\Omega)},
        \end{align*}
        which holds for any $\phi \in H^1(\Omega)$, where $C_P>0$ denotes the Poincaré constant on $\Omega$.
        Applying this inequality with $\phi = (r-\ac)(s)$, we obtain
        \begin{align}\label{lemma : generalized Poincaré - singular limit : III}
            \|(r-\ac)(s)\|_{L^2(\Omega)}
            \leq 
            C_P \|\nablax r(s) - \nablax \ac(s)\|_{L^2(\Omega)}
            +
            \Bigl\|\fint_\Omega (r-\ac)(s)\Bigl\|_{L^2(\Omega)}.
        \end{align}
        With relation $\eqref{lemma : generalized Poincaré - singular limit : II}$ and Jensen's inequality, we estimate the second term on the right-hand side of inequality $\eqref{lemma : generalized Poincaré - singular limit : III}$ as
        \begin{align}\label{lemma : generalized Poincaré - singular limit : IV}
            \Bigl\|\fint_\Omega (r-\ac)(s)\Bigl\|_{L^2(\Omega)}
            \leq
            \frac{\paracoup}{\coup}
            \Bigl\|\fint_\Omega \partial_t \ac(s) \Bigl\|_{L^2(\Omega)}
            +
            ||e_\coup||_{L^2(\Omega)}
            \leq
            \frac{\paracoup}{\coup}
            \|\partial_t \ac(s)\|_{L^2(\Omega)}
            +
            |e_\coup| \sqrt{|\Omega|}.
        \end{align}
        Combining $\eqref{lemma : generalized Poincaré - singular limit : III}$ and $\eqref{lemma : generalized Poincaré - singular limit : IV}$ yields $\eqref{lemma : generalized Poincaré - singular limit : inequality}$.
    \end{proof}
    We have now everything prepared in order to prove our second main result Theorem~\ref{thm main result:singular limit}.
    \begin{proof}[Proof of Theorem~\ref{thm main result:singular limit}]
        In order to ease the notation, we set
        \begin{align*}
            \curlyE_\coup
            :=
            \curlyE(\arho,\auvec,\ac\mid r,\Uvec),
            \quad
            \curlyE_\coup^0
            :=
            \curlyE(\arho^0,\auvec^0,\ac^0\mid r^0,\Uvec^0),
        \end{align*}
        where the relative energy $\curlyE(\arho,\auvec,\ac \mid r,\Uvec)$ is defined as in $\eqref{defi Erel for relax limit}$.
        By Proposition~\ref{proposition : REI singular limit} we have that the relative energy inequality
        \begin{align}
            &\curlyE_\coup(s)
            +
            \int_0^s \int_\Omega 
            \lr{
            \Svisc(\nablax\auvec-\nablax\Uvec):(\nablax\auvec-\nablax\Uvec)
            +
            \paracoup |\partial_t\ac|^2 
            }
            \, \dd x \, \dd t\nonumber
            \\
            &\qquad \leq 
            \curlyE_\coup^0
            +
            \mathcal{R}_\coup^\mathrm{NSE}
            +
            \sum\limits_{i=1}^7
            \mathcal{J}_\coup^i\label{singular limit Pf REI}
        \end{align}
        holds for almost all $s \in (0,T)$, where $\mathcal{R}_\coup^{\mathrm{NSE}}$ and $\mathcal{J}_\coup^i$ for $i \in \{1,2,\cdots,7\}$ are defined as in Proposition~\ref{proposition : REI singular limit}.
        Our goal is to estimate the right-hand side of this inequality appropriately, so that we can apply Gronwall's inequality.
        For $\mathcal{R}^{\mathrm{NSE}}_\coup$ we perform the same estimates as in the proof of Theorem~\ref{thm main result:W-S uniqueness} (cf. $\eqref{Pf WS-uniqueness III}$) and obtain 
        \begin{align}\label{Pf singular limit RNSE}
            |\mathcal{R}^{\mathrm{NSE}}_\coup|
            \leq 
            K_1 \int_0^s 
            \curlyE_\coup(t)
            \, \dd t
            +
            \frac{1}{4}
            \int_0^s \int_\Omega 
            \Svisc(\nablax\uvec-\nablax\Uvec) : (\nablax \uvec - \nablax\Uvec)
            \, \dd x\, \dd t,
        \end{align}
        for some positive constant $K_1$ that does not depend on $\coup$.
        For $\mathcal{J}^1_\coup$ we first infer by Young's inequality
        \begin{align}
            \int_0^s \int_\Omega 
            |\arho-r|^2 
            \, \dd x \, \dd t
            &\leq 
            \int_0^s \int_\Omega 
            \lr{
            2|\arho - \ac |^2
            +
            2|\ac-r|^2
            }
            \, \dd x \, \dd t\nonumber
            \\
            &\leq 
            \frac{4}{\coup}
            \int_0^s
            \mathcal{E}_\coup(t)
            \, \dd t
            +
            2
            \int_0^s\int_\Omega 
            |\ac-r|^2 
            \, \dd x \, \dd t.\label{Pf singular limit J^1 I}
        \end{align}
        By inequality $\eqref{lemma : generalized Poincaré - singular limit : inequality}$ in Lemma~\ref{lemma : generalized Poincaré - singular limit}, we estimate
        \begin{align}
            &\int_0^s \int_\Omega 
            |\ac-r|^2 
            \, \dd x \, \dd t\nonumber\\
            &\leq 
            2C_P^2
            \int_0^s \int_\Omega 
            |\nablax\ac - \nablax r|^2 
            \, \dd x \, \dd t
            +
            \frac{2\paracoup^2}{\coup^2}
            \int_0^s \int_\Omega 
            |\partial_t\ac|^2
            \, \dd x \, \dd t
            +
            2
            \int_0^s
            e_\coup^2 |\Omega|
            \, \dd t\nonumber
            \\
            &\leq 
            \frac{4C_P^2}{\kappa} 
            \int_0^s
            \mathcal{E}(s)
            \, \dd t
            +
            \frac{2\paracoup}{\coup^2}
            \int_0^s \int_\Omega 
            \paracoup|\partial_t \ac |^2 
            \, \dd x \, \dd t
            +
            2e_\coup^2 |\Omega| T.\label{Pf singular limit J^1 II}
        \end{align}
        With $\eqref{Pf singular limit J^1 I}$, $\eqref{Pf singular limit J^1 II}$, the Lipschitz continuity of $q$ and Young's inequality, we estimate $\mathcal{J}^1_\coup$ for any $\coup\geq 1$ as
        \begin{align}
            |\mathcal{J}^1_\coup|
            &\leq 
            \frac{L_q^2}{\lambda}
            \int_0^s \int_\Omega 
            |\arho-r|^2 
            \, \dd x \, \dd t
            +
            \frac{\lambda}{4}
            \int_0^s \int_\Omega 
            |\div\uvec-\div\Uvec|^2 
            \, \dd x \, \dd t\nonumber
            \\
            &\leq 
            K_2
            \int_0^s 
            \mathcal{E}(t)
            \, \dd t
            +
            k_\coup^1
            \int_0^s \int_\Omega 
            \paracoup|\partial_t\ac|^2 
            \, \dd x \, \dd t
            +
            l_\coup^1\nonumber\\
            &\qquad+
            \frac{1}{4}
            \int_0^s \int_\Omega 
            \Svisc(\nablax\uvec - \nablax\Uvec) : (\nablax\uvec - \nablax\Uvec)
            \, \dd x \, \dd t,\label{Pf singular limit J^1}
        \end{align}
        with
        \begin{align*}
            K_2
            :=
            \frac{4L_q^2}{\lambda}
            \biggl(1+\frac{2C_P^2}{\kappa}\biggr),
            \quad 
            k_\coup^1
            :=
            \frac{4L_q^2\paracoup}{\lambda\coup^2},
            \quad
            l_\coup^1
            :=
            \frac{4L_q^2|\Omega|e_\coup^2T}{\lambda}.
        \end{align*}
        For $\mathcal{J}^2_\coup$, we use integration by parts, Young's inequality and inequality $\eqref{Pf singular limit J^1 II}$ to estimate
        \begin{align}
            | 
            \mathcal{J}_\coup^2
            |
            &\leq 
            \int_0^s \int_\Omega
            \lr{
            \frac{\kappa}{2} |\div\Uvec|\,|\nablax(r-\ac)|^2 
            +
            \kappa |r-\ac|\,|\nablax(r-\ac)\cdot\nablax\div
            \Uvec|
            }
            \, \dd x\, \dd t\nonumber
            \\
            &\leq 
            \Bigl(\|\div\Uvec\|_{\curlyC(\closOmT)} + \|\nablax\div\Uvec\|_{\curlyC(\closOmT)}^2 \Bigr)
            \int_0^s \curlyE_\coup(t)
            \, \dd t
            +
            \int_0^s \int_\Omega
            \frac{\kappa}{2}|r-\ac|^2 
            \, \dd x \, \dd t\nonumber
            \\
            &\leq 
            \Bigl( \|\div\Uvec\|_{\curlyC(\closOmT)} + \|\nablax\div \Uvec\|_{\curlyC(\closOmT)}^2 + 2 C_P^2\Bigr) 
            \int_0^s
            \curlyE_\coup(t) 
            \, \dd t\nonumber
            \\
            &\qquad+
            \frac{\kappa\paracoup}{\coup^2}
            \int_0^s \int_\Omega 
            \paracoup
             |\partial_t\ac|^2
            \, \dd x \, \dd t
            +
            \kappa|\Omega| e_\coup^2 T\nonumber
            \\
            &\leq 
            K_3 
            \int_0^s 
            \curlyE(t)
            \, \dd t
            +
            k_\coup^2 
            \int_0^s \int_\Omega 
            \paracoup|\partial_t \ac|^2 
            \, \dd x \, \dd t
            +
            l_\coup^2\label{Pf singular limit J^2}
        \end{align}
        with
        \begin{align*}
            K_3
            :=
            \Bigl( 
            \|\div\Uvec\|_{\curlyC(\closOmT)}
            +
            \|\nablax \div\Uvec\|_{\curlyC(\closOmT)}^2
            +
            2C_P^2 
            \Bigr),
            \quad
            k_\coup^2 
            :=
            \frac{\kappa\paracoup}{\coup^2},
            \quad
            l_\coup^2 := \kappa|\Omega| e_\coup^2 T.
        \end{align*}
        For $\mathcal{J}_\coup^3$ and $\mathcal{J}^5_\coup$, we estimate
        \begin{align}
            |
            \mathcal{J}^3_\coup 
            |
            +
            |
            \mathcal{J}^5_\coup
            |
            \leq 
            K_4
            \int_0^s 
            \curlyE_\coup(t)
            \, \dd t\label{Pf singular limit J^3 and J^5}
        \end{align}
        with
        \begin{align*}
            K_4 := 6\|\nablax\Uvec\|_{C(\closOmT)} + \|\div\Uvec\|_{\curlyC(\closOmT)}.
        \end{align*}
        For $\mathcal{J}_\coup^4$, we estimate with Young's inequality
        \begin{align}
            |\mathcal{J}_\coup^4|
            &\leq 
            \int_0^s \int_\Omega 
            \Biggl(
            \frac{\kappa ^2\coup}{2} |\arho-\ac|^2 
            +
            \frac{\|\Uvec\|_{\curlyC(\closOmT)}^2 \|\nablax\Deltax r\|_{\curlyC(\closOmT)}^2}{2\coup}
            \Biggr)
            \, \dd x \, \dd t\nonumber
            \\
            &\leq 
            K_5 
            \int_0^s
            \curlyE(t)
            \, \dd t
            +
            l_{\coup}^3,\label{Pf singular limit J^4}
        \end{align}
        with
        \begin{align*}
            K_5
            := 
            \kappa^2, 
            \quad 
            l_3^\coup
            :=
            \frac{\|\Uvec\|_{\curlyC(\closOmT)}^2 \|\nablax\Deltax r\|_{\curlyC(\closOmT)}^2 |\Omega| T}{2\coup}.
        \end{align*}
        For $\mathcal{J}^6_\coup$, we estimate with Young's inequality and inequality $\eqref{Pf singular limit J^1 II}$
        \begin{align}
            |\mathcal{J}_\coup^6|
            &\leq 
            \int_0^s \int_\Omega 
            |\paracoup\partial_t\ac \div(\ac \Uvec)| 
            \, \dd x \, \dd t\nonumber
            \\
            &\leq 
            \int_0^s \int_\Omega 
            \lr{
            |\paracoup\partial_t\ac \div\bigl((\ac -r) \Uvec\bigr)| 
            +
            |\paracoup\partial_t \ac \div(r\Uvec)|
            }
            \, \dd x \, \dd t\nonumber
            \\
            &\leq 
            \int_0^s \int_\Omega 
            \lr{
            \frac{\paracoup^2}{2}|\partial_t\ac |^2 
            +
            \frac{1}{2}|\div\bigl( (\ac-r)\Uvec|^2
            +
            \frac{\paracoup}{4}|\partial_t\ac|^2 
            +
            \paracoup\|\div(r\Uvec)\|_{\curlyC(\closOmT)}^2
            }
            \, \dd x \, \dd t\nonumber
            \\
            &\leq 
            \int_0^s \int_\Omega
            \lr{
            \Bigl(
            \frac{\paracoup}{2}
            +
            \frac{1}{4}\Bigr)
            \paracoup
            |\partial_t\ac |^2
            +
            \frac{1}{2}
            \|\Uvec\|_{\curlyC(\closOmT)}^2 |\nablax (\ac - r)|^2 
            }
            \, \dd x \, \dd t\nonumber\\
            &\qquad+
            \int_0^s \int_\Omega
            \frac{1}{2}
            \|\div\Uvec\|_{\curlyC(\closOmT)}^2 |\ac - r|^2 
            \, \dd x \, \dd t
            +
            \paracoup \|\div(r\Uvec)\|_{\curlyC(\closOmT)}^2 
            |\Omega| T\nonumber
            \\
            &\leq
            K_6
            \int_0^s 
            \curlyE(t)
            \, \dd t
            +
            k_\coup^3 \int_0^s 
            \int_\Omega 
            \paracoup|\partial_t \ac|^2 
            \, \dd x
            \, \dd t
            +
            l_\coup^4\label{Pf singular limit J^6}
        \end{align}
        with
        \begin{align*}
            &K_6 
            := 
            \frac{ \|\Uvec\|_{\curlyC(\closOmT)}^2 + 2C_P^2 \|\div\Uvec\|_{\curlyC(\closOmT)}^2}{\kappa},
            \quad 
            k_\coup^3
            :=
            \frac{1}{4}
            +
            \frac{\paracoup}{2}
            +
            \frac{\paracoup}{\coup^2}\|\div \Uvec\|_{\curlyC(\closOmT)}^2,
            \\ 
            &l_\coup^4
            :=
            \Bigl( \|\div\Uvec\|_{\curlyC(\closOmT)}^2e_\coup^2 
            +
            \paracoup\|\div(r\Uvec)\|_{\curlyC(\closOmT)}^2 \Bigr)|\Omega|T.
        \end{align*}
        For $\mathcal{J}_\coup^7$, we first rewrite by using the continuity equation $\eqref{def:finite energy weak sol:continuity}$ and the regularity of $\ac$ 
        \begin{align*}
            \mathcal{J}_\coup^7
            &=
            \Biggl[ \int_\Omega   \ac \kappa \Deltax r \, \dd x \Biggr]_{\tau=0}^{\tau=s}
            -
            \int_0^s \int_\Omega 
            \kappa \partial_t \Deltax  r \, \ac 
            \, \dd x \, \dd t
            -
            \Biggl[ \int_\Omega \arho \kappa \Deltax r \, \dd x \Biggr]_{\tau=0}^{\tau=s}\\
            &\qquad+
            \int_0^s \int_\Omega 
            \kappa \partial_t \Deltax r\, \arho
            \, \dd x \, \dd t
            \\
            &=
            \int_\Omega
            \lr{
            (\arho^0-\ac^0)\kappa \Deltax r^0
            -
            (\arho-\ac)(s) \kappa \Deltax r(s) 
            }
            \, \dd x
            +
            \int_0^s \int_\Omega 
            (\arho - \ac) \kappa \partial_t\Deltax r
            \, \dd x \, \dd t.
        \end{align*}
        Then, we estimate with Young's inequality
        \begin{align}
            |\mathcal{J}_\coup^7|
            &\leq 
            \int_\Omega 
            \lr{
            \frac{\coup}{2} |\arho^0-\ac^0|^2 
            +
            \frac{\kappa^2}{2\coup}|\Deltax r^0|^2 
            +
            \frac{\coup}{4}|\arho(s)-\ac(s)|^2
            +
            \frac{\kappa^2}{\coup}|\Deltax r(s)|^2 
            }
            \, \dd x\nonumber
            \\
            &\qquad
            +
            \int_0^s \int_\Omega 
            \lr{
            \frac{\coup}{2} |\arho-\ac|^2 
            +
            \frac{\kappa^2}{2\coup} |\partial_t\Deltax r|^2 
            }
            \, \dd x \, \dd t\nonumber
            \\
            &\leq 
            \int_0^s
            \curlyE_\coup(t) 
            \, \dd t
            +
            \frac{1}{2} \curlyE_\coup(s)
            +
            l_\coup^5,\label{Pf singular limit J^7}
        \end{align}
        using
        \begin{align*}
            l_\coup^5
            :=
            \curlyE_\coup^0 
            +
            \frac{\kappa^2 |\Omega| T}{2\coup}\|\partial_t\Deltax r\|_{\curlyC(\closOmT)}^2 
            +
            \frac{3\kappa^2|\Omega|}{2\coup}\|\Deltax r\|^2_{\curlyC(\closOmT)}.
        \end{align*}
        Combining $\eqref{singular limit Pf REI}, \eqref{Pf singular limit RNSE}, \eqref{Pf singular limit J^1}, \eqref{Pf singular limit J^2}, \eqref{Pf singular limit J^3 and J^5}, \eqref{Pf singular limit J^4}$, $\eqref{Pf singular limit J^6}$ and $\eqref{Pf singular limit J^7}$ yields for any $\coup \geq 1$ that
        \begin{align}
            &\frac{1}{2}\curlyE_\coup(s)
            +
            \frac{1}{2}
            \int_0^s \int_\Omega 
            \Svisc(\nablax\auvec-\nablax\Uvec):(\nablax\auvec-\nablax\Uvec)
            \, \dd x \, \dd t
            +
            (1-k_\coup^\Sigma)\int_0^s \int_\Omega \paracoup|\partial_t\ac|^2 
            \, \dd x \, \dd t\nonumber
            \\
            &\leq 
            K_\Sigma
            \int_0^s 
            \curlyE_\coup(t)
            \, \dd t
            +
            l_\coup^\Sigma \label{singular limit Pf prep Gronwall}
        \end{align}
        with $K_\Sigma := 1 + \sum\limits_{i=1}^{6}K_i$ being a constant that does not depend on $\coup$, $l_\coup^\Sigma :=  \curlyE_\coup^0 +\sum\limits_{i=1}^{5}l_\coup^i$, and
        \begin{align*}
            k_\coup^\Sigma 
            :=
            (k_\coup^1+k_\coup^2+k_\coup^3)
            =
            \frac{1}{4}
            +
            \paracoup\biggl(
            \frac{4L_q^2}{\lambda\coup^2}
            +
            \frac{\kappa+\|\div\Uvec\|_{\curlyC(\closOmT)}^2}{\coup^2}
            +
            \frac{1}{2}
            \biggr).
        \end{align*}
        Since $\paracoup=\paracoup(\coup) \to 0$ for $\coup \to \infty$, there exists some $\coup_0\geq1$, such that for any $\coup \geq \coup_0$, we have that $k_\coup^\Sigma \leq \frac{1}{2}$.
        Thus, inequality $\eqref{singular limit Pf prep Gronwall}$ implies for any $\coup\geq \coup_0$ that the inequality
        \begin{align}
            &\curlyE_\coup(s)
            +
            \int_0^s \int_\Omega 
            \Svisc(\nablax\auvec-\nablax\Uvec):(\nablax\auvec-\nablax\Uvec)
            \, \dd x \, \dd t
            +
            \int_0^s \int_\Omega \paracoup|\partial_t\ac|^2 
            \, \dd x \, \dd t\nonumber
            \\
            &\leq 
            2K_\Sigma
            \int_0^s 
            \curlyE_\coup(t)
            \, \dd t
            +
            2l_\coup^\Sigma\label{singular limit Pf FINAL 0}
        \end{align}
        holds for almost all $s \in (0,T)$.
        Applying Gronwall's inequality yields then for $\coup\geq \coup_0$ that the inequality
        \begin{align}\label{singular limit Pf FINAL I}
            \curlyE_\coup(s)
            \leq 
            2l_\coup^\Sigma 
            \Bigl( 
            1 + 2K_\Sigma \exp(2K_\Sigma T)T 
            \Bigr) 
        \end{align}
        holds for almost all $s\in(0,T)$.\\
        To ease the notation, we denote from now on by $\cvconst>0$ a generic positive constant that may vary from line to line, but does only depend on $\mu,\lambda,\kappa,\gamma,|\Omega|,T,r_{\mathrm{min}}$ and on the norms in $\eqref{thm main result:singular limit:dependencies}$.
        With this convention, we have by definition of $l_\coup^\Sigma$ and $K^\Sigma$ that
        \begin{align*}
            l_\coup^\Sigma \leq \cvconst\cvrate,
            \quad K_\Sigma \leq \cvconst,
        \end{align*}
        where $\cvrate$ is defined as in $\eqref{definition cvrate and ealpha}$, and thus, by $\eqref{singular limit Pf FINAL I}$,
        \begin{align}\label{singular limit Pf FINAL estimate E_coup}
            \|\mathcal{E}_\coup\|_{L^\infty(0,T)} \leq \cvconst\cvrate.
        \end{align}
        With Lemma~\ref{lemma:auxiliary convexity estimates}, we conclude from $\eqref{singular limit Pf FINAL estimate E_coup}$ that
        \begin{align}
            &\|\sqrt{\arho}(\auvec-\Uvec)\|_{L^\infty(0,T;L^2(\Omega;\RR^3))}^2
            +
            \|\arho - r\|_{L^\infty(0,T;L^{\hat{\gamma}}(\Omega))}^2
            +
            \coup \|\arho - \ac \|_{L^\infty(0,T;L^2(\Omega))}^2\nonumber
            \\
            &\quad+
            \|\nablax \ac - \nablax r\|_{L^\infty(0,T;L^2(\Omega;\RR^3))}^2 \leq \cvconst \cvrate,\label{singular limit Pf FINAL estimte from erel}
        \end{align}
        where $\hat{\gamma}:=\min\{2,\gamma\}$.
        In particular, we have by Hölder's inequality for almost all $s \in (0,T)$ that
        \begin{align*}
            \biggl|\int_\Omega (\ac(s)-r(s))\, \dd x \biggr|
            &\leq 
            \int_\Omega |\ac(s)-\arho(s)|\, \dd x 
            +
            \int_\Omega |\arho(s) - r(s)|\, \dd x
            \\
            &
            \leq 
            \cvconst\|\ac(s)-\arho(s)\|_{L^2(\Omega)}
            +
            \cvconst \|\arho(s) - r(s)\|_{L^{\hat{\gamma}}(\Omega)}
            \\
            &\leq \cvconst \sqrt{\cvrate}.
        \end{align*}
        With $\eqref{singular limit Pf FINAL estimte from erel}$ and Poincaré's inequality we deduce from this estimate
        \begin{align}
            \|\ac-r\|_{L^\infty(0,T;L^2(\Omega))}
            &\leq \cvconst\|\nablax \ac - \nablax r\|_{L^\infty(0,T;L^2(\Omega;\RR^3))}
            +
            \cvconst\Biggl\| \int_\Omega (\ac-r)\, \dd x \Biggr\|_{L^\infty((0,T))}
            \nonumber
            \\
            &\leq 
            \cvconst \sqrt{\cvrate}.\label{singular limit Pf FINAL estimate ac - r}
        \end{align}
        Combining $\eqref{singular limit Pf FINAL estimte from erel}$ and $\eqref{singular limit Pf FINAL estimate ac - r}$ yields
        \begin{align*}
            \|\arho - r\|^2_{L^\infty(0,T;L^2(\Omega))}
            \leq 
            \cvconst\|\arho-\ac\|_{L^\infty(0,T;L^2(\Omega))}^2 
            +
            \cvconst\|\ac - r\|_{L^\infty(0,T;L^2(\Omega))}^2
            \leq 
            \cvconst
            \cvrate.
        \end{align*}
        Going back to $\eqref{singular limit Pf FINAL 0}$, we obtain by using $\eqref{singular limit Pf FINAL estimate E_coup}$, Korn's and Poincaré's inequality
        \begin{align*}
            \|\auvec-\Uvec\|_{L^2(0,T;H^1(\Omega;\RR^3))}^2
            +
            \paracoup(\coup) \|\partial_t\ac\|_{L^2(0,T;L^2(\Omega))}^2
            \leq 
            \cvconst\cvrate,
        \end{align*}
        which completes the proof of $\eqref{thm main result:singular limit:convergence rates}$.\\
        If we assume for the initial data the relation $\eqref{thm main result:singular limit:assumption initial data}$, then we have by Lemma~\ref{lemma:auxiliary convexity estimates} in particular
        \begin{align*}
            \lim\limits_{\coup \to \infty}|e_\coup| =0.
        \end{align*}
        Thus, we have
        \begin{align*}
            \lim\limits_{\coup\to\infty} \cvrate = 0
        \end{align*}
        and the convergences $\eqref{thm main result:singular limit:convergences}$ follow from $\eqref{thm main result:singular limit:convergence rates}$, since $\cvconst$ does not depend on $\coup$.
        The proof is now complete.
    \end{proof}
    \section{Conclusions}\label{Sec:Conclusions}
    In this paper we considered a relaxation formulation for the compressible NSK equations with a pressure function of Van-der-Waals type that was proposed in \cite{HitzKeimMunzRohde2020} as a phase transition model for an isothermal viscous compressible two-phase flow.    
    We investigated the model from a rigorous mathematical perspective on a bounded domain and verified two important analytical properties.
    As a first result, we showed that finite energy weak solutions for the initial-boundary value problem of the relaxed NSK equations satisfy the weak--strong uniqueness property.
    Motivated by the great importance for the numerical analysis of the governing equations, it seems very interesting to find the dissipative measure-valued formulation of the relaxed NSK equations and to prove the weak-strong uniqueness property for this weaker class of solutions.
    With such a result, a rigorous numerical analysis as performed in \cite{FeireislBook2021} would be possible.
    As our second result, we proved that the relaxed NSK equations approach the NSK equations in the relaxation limit $\coup \to \infty$ and $\paracoup\to 0$ in certain strong norms, provided a classical solution of the NSK equations exist.
    As a byproduct we identified a convergence rate for the corresponding norms.
    With this result, we rigorously justify the relaxed NSK equations as an approximate model for the NSK equations complementing the numerical experiments and the formal argument given in \cite{HitzKeimMunzRohde2020}.
    In both of our result, we postulated the global-in-time existence of finite energy weak solutions to the initial-boundary value problem $\eqref{relaxed NSK field eq continuity}$--$\eqref{relaxed NSK IC}$ emanating from initial data that satisfy the regularity $\eqref{defi FEW initial conditions}$.
    To complement our results, a rigorous realization of such a global-in-time existence result is planned as forthcoming work.
    \par
    Both of our results apply for pressure functions satisfying $\eqref{assump on p:decomposition}$--$\eqref{assump on p:monotonicity h}$ and thus, for pressure functions of Van-der-Waals type.
    In particular, we allow for a two-phase setting in both results.
    However, we require that the pressure function is defined on the whole interval $[0,\infty)$ and thus, our theory does not account for the original Van-der-Waals pressure function, that is,
    \begin{align*}
        p_{\mathrm{VdW}}\colon [0,b) \to [0,\infty),
        \quad \rho \mapsto \frac{RT_{\mathrm{ref}}\rho}{b-\rho} - a \rho^2,
    \end{align*}
    where $T_\mathrm{ref} >0$ denotes the reference temperature, $R>0$ denotes the specific gas constant and $a,b>0$ denote material coefficients.
    The Van-der-Waals pressure function is singular for $\rho \to b$ and belongs to the class of hard-sphere pressure functions.
    For certain pressure functions in this class, the weak-strong uniqueness principle for the compressible Navier--Stokes equations was verified by the first author in \cite{Chaudhuri2019}.
    It would be very interesting to generalize our results to the class of hard-sphere pressure functions.
    Both of our results rely strongly on the relative energy inequality derived in Section~\ref{Sec:Relative Energy Inequality}, which is of independent interest to perform different singular limits.
    It would be highly interesting to use this relative energy inequality to perform other physically reasonable singular limits as the low Mach number limit or the sharp interface limit for the relaxed NSK equations. 
    \section*{Acknowledgements}
    Funded by Deutsche Forschungsgemeinschaft (DFG, German Research Foundation) under Germany´s Excellence Strategy – EXC 2075 – 390740016. 
    The research of N.C. is supported by the “Excellence Initiative Research University (IDUB)” program at the University of Warsaw. 
    The authors acknowledge support from the Institute of Mathematics of the Czech Academy of Sciences. They thank Eduard Feireisl for the hospitality and the valuable suggestions.
    FW thanks Mária Luká\v{c}ová-Medvid'ová for the hospitality at the JGU Mainz and the inspiring discussions.
\bibliographystyle{plain}
\bibliography{sample}
\end{document}